\documentclass[10pt,journal,compsoc, twoside]{IEEEtran}
\ifCLASSOPTIONcompsoc
  \usepackage[nocompress]{cite}
\else
  \usepackage{cite}
\fi
\usepackage{xcolor}
\ifCLASSINFOpdf
  \usepackage[pdftex]{graphicx}
\else
  \usepackage[dvips]{graphicx}
\fi
\usepackage{pgfplots} 
\pgfplotsset{width=7cm,compat=1.14}
\usepgfplotslibrary{statistics}
\usepackage{ifthen}
\usepackage{tikz} 
\usetikzlibrary{calc,patterns,shapes}
\usetikzlibrary{arrows.meta}
\tikzset{>={Latex[width=2mm,length=2mm]}}
\usetikzlibrary{decorations.pathreplacing}
\usepgfplotslibrary{fillbetween}
\usetikzlibrary{patterns}
\usepackage{booktabs}
\usepackage{siunitx}
\usepackage{verbatim}

\usepackage{microtype} 
\usepackage{printlen}
\usepackage{layouts}
\usepackage[normalem]{ulem}
\usepackage{amsmath}
\usepackage{algorithmic}
  \usepackage[caption=false,font=footnotesize,labelfont=sf,textfont=sf]{subfig}
\usepackage{url}

\usepackage[hidelinks]{hyperref}

\newcommand{\bq}{\begin{equation}}
\newcommand{\eq}{\end{equation}}
\newcommand{\bytes}{\mbox{bytes}}
\newcommand{\byte}{\mbox{byte}}
\newcommand{\second}{\mbox{s}}

\newcommand{\flop}{\mbox{flop}}

\newcommand{\bits}{\mbox{bits}}

\newcommand{\GBS}{\mbox{G\byte/\second}}

\newcommand{\GFS}{\mbox{G\flop/\second}}

\newcommand{\BF}{\mbox{\byte/\flop}}

\newcommand{\MB}{\mbox{MB}}

\newcommand{\MiB}{\mbox{MiB}}
\newcommand{\KiB}{\mbox{KiB}}

\newcommand{\eos}{\;.}
\newcommand{\cma}{\;,}
\newcommand{\rlm}{roof{}line model}
\newcommand{\rl}{roof{}line}

\newcommand{\NNZR}{\mbox{$N_\mathrm{nzr}$}}
\newcommand{\NR}{\mbox{$N_\mathrm{r}$}}
\newcommand{\NNZ}{\mbox{$N_\mathrm{nz}$}}

\newcommand{\likwid}{\texttt{LIKWID}}
\newcommand{\likwidperfctr}{\texttt{likwid-perfctr}}
\newcommand{\likwidpin}{\texttt{likwid-pin}}
\newcommand{\likwidbench}{\texttt{likwid-bench}}

\definecolor{tumbleweed}{rgb}{0.87, 0.67, 0.53}

\newcommand{\spmv}{SpMV}

\newcommand{\MPK}{\textsc{MPK}}

\newcommand*{\rom}[1]{\expandafter\@slowromancap\romannumeral #1@}

\newenvironment{customlegend}[1][]{%
	\begingroup
	\csname pgfplots@init@cleared@structures\endcsname
	\pgfplotsset{#1}%
}{%
	\csname pgfplots@createlegend\endcsname
	\endgroup
}%
\def\addlegendimage{\csname pgfplots@addlegendimage\endcsname}

\definecolor{applegreen}{rgb}{0.55, 0.71, 0.0}
\definecolor{amethyst}{rgb}{0.6, 0.4, 0.8}
\definecolor{amber}{rgb}{1.0, 0.75, 0.0}

\definecolor{col0}{rgb}{  0.5508,    0.8242,    0.7773 }
\definecolor{col1}{rgb}{  0.98,    0.93,    0.36 }
\definecolor{col2}{rgb}{  0.7422,    0.7266,    0.8516 }
\definecolor{col3}{rgb}{  0.9805,    0.5000,    0.4453 }
\definecolor{col4}{rgb}{  0.5000,    0.6914,    0.8242 }
\definecolor{col5}{rgb}{  0.9883,    0.7031,    0.3828 }
\definecolor{col6}{rgb}{  0.6992,    0.8672,    0.4102 }
\definecolor{col7}{rgb}{  0.9844,    0.8008,    0.8945 }
\definecolor{col8}{rgb}{  0.6484,    0.8047,    0.8867 }
\definecolor{col9}{rgb}{  0.1211,    0.4688,    0.7031 }
\definecolor{col10}{rgb}{  0.6953,    0.8711,    0.5391 }
\definecolor{col11}{rgb}{  0.1992,    0.6250,    0.1719 }
\definecolor{col12}{rgb}{  0.9805,    0.6016,    0.5977 }
\definecolor{col13}{rgb}{  0.8867,    0.1016,    0.1094 }
\definecolor{col14}{rgb}{  0.9883,    0.7461,    0.4336 }

\begin{document}
%
\title{Level-based Blocking for Sparse Matrices: Sparse Matrix-Power-Vector Multiplication}

\author{Christie~Alappat,
	Georg~Hager,
	Olaf~Schenk
	and~Gerhard~Wellein
	\IEEEcompsocitemizethanks{\IEEEcompsocthanksitem C. Alappat, G. Hager and G. Wellein are with Erlangen National High Performance Computing Center at Friedrich-Alexander-Universität Erlangen-Nürnberg.
		E-mail: \{christie.alappat, georg.hager, gerhard.wellein\}@fau.de
		\IEEEcompsocthanksitem G. Wellein is also with the Department of Computer Science, Friedrich-Alexander-Universität Erlangen-Nürnberg.
		\IEEEcompsocthanksitem O. Schenk is with the Institute of Computing at Faculty of Informatics, Universit\`{a} della Svizzera italiana.
E-mail: {olaf.schenk@usi.ch}	
}
}	

%
%

\markboth{}{\protect{Alappat \MakeLowercase{\textit{et al.}}: Level-based blocking for sparse matrices}}
%



\IEEEtitleabstractindextext{%
\begin{abstract}

The multiplication of a sparse matrix with a dense vector (SpMV) is a key component in many numerical schemes and its performance is known to be severely limited by main memory access. 
Several numerical schemes require the multiplication 
of a sparse matrix polynomial with a dense vector which is typically implemented as a sequence of SpMVs.
This results in low performance and ignores the potential to increase the arithmetic intensity by reusing the matrix data from cache.  
In this work we use the recursive algebraic coloring engine (RACE) to enable blocking of sparse matrix data across the polynomial computations. 
In the graph representing the sparse matrix we form levels using a breadth-first search. 
Locality relations of these levels are then used to improve  spatial and temporal locality when accessing the matrix data and to implement an efficient multithreaded parallelization.
Our approach is independent of the matrix structure and avoids shortcomings of
existing ``blocking'' strategies in terms of hardware efficiency and parallelization overhead.
We quantify the quality of our implementation using performance modelling and demonstrate speedups of up to 3$\times$ and 5$\times$ compared to an optimal SpMV-based baseline 
on a single multicore chip of recent Intel and AMD architectures. 
As a potential application, we demonstrate the benefit of our implementation for a Chebyshev time propagation scheme, representing the class of polynomial approximations to exponential integrators. 
Further numerical schemes which may benefit from our developments include $s$-step Krylov solvers and power clustering algorithms.

\end{abstract}

}

\maketitle

\IEEEdisplaynontitleabstractindextext

%
\IEEEpeerreviewmaketitle

\IEEEraisesectionheading{\section{Introduction and Related Work}\label{sec:introduction}}

Sparse matrix-vector multiplication (SpMV) is a critical building block for a wide variety of computational algorithms used in science, engineering, and data analytics. 
The SpMV kernel is known to perform poorly on modern compute devices due to its low arithmetic intensity and often irregular memory access pattern.
Most performance optimization efforts target a single SpMV invocation. 
To minimize the data access costs to the matrix entries, a plethora of data layout choices have been proposed for GPGPUs~\cite{GPU_review} 
and CPUs~\cite{barrett_formats,Vuduc2003:thesis,CSB},
including hardware-agnostic formats~\cite{KreutzerSELLC14}.
These formats typically ensure linear access to matrix data, but the input vector is always accessed indirectly and therefore potentially in an irregular way.
Optimization strategies like matrix reordering or partitioning  techniques~\cite{SpMV_reordering} aim to reduce the reuse distances in the vector accesses and thus improve the performance. 
Finally, at the kernel implementation level, automatic performance optimization 
for SpMV has been a subject of research for decades. 
These approaches mainly account for the complexity of cache-based microprocessors, where SpMV performance maybe extremely sensitive to the spatial/temporal data access locality, out-of-order instruction capability, register scheduling, and SIMD vectorization. 
Choosing parameters for these code optimizations and choosing among alternative implementations is critical for efficient hardware utilization.
It has been demonstrated~\cite{BEBOP2006,Vuduc2003:thesis, balaprakash2018autotuning, hong2019adaptive} 
that it is possible to build an automatic tuning system capable of generating implementations that are on par with or even outperform the best manually tuned code.

In this work, we extend SpMV performance tuning research towards automatic data reuse optimization across several SpMV invocations in the \emph{sparse matrix-power-vector kernel} (MPK), 
which computes $Ax$, $A^2x, A^3x, \cdots, A^kx$ for matrix $A$, vector $x$, 
and a small constant $k$.
Our focus is on thread-level parallel and efficient CPU implementation of MPK 
using 
the popular compressed row storage (CRS)
sparse matrix format.
To this end we extend the recursive algebraic coloring engine (RACE) 
framework~\cite{RACE} to tackle the dependencies between several SpMV invocations
in the MPK.
The algebraic formulation used in RACE is general in the sense that it does not
assume any special structure in the underlying matrix.

The need for software 
implementations and structures for MPK is exemplified by communication-avoiding algorithms~\cite{demmel_MPK_tr,hoemmen_thesis,Carson:EECS-2015-179,7914608},
which have been proposed to improve performance by trading redundant computation for memory traffic.  
In these algorithms, independent SpMV invocations are replaced by 
the MPK to compute $A^kx$.
Once the computation has been performed, the next $k$ steps of the solver can proceed without further memory accesses to $A$ by combining vectors from this set. 

There has been some research in exploiting data locality in MPK, mostly motivated by classic blocking strategies well established in stencil computations.
In particular, in~\cite{mohiyudeen}
blocking schemes for MPK have been developed that first partition 
the graph of a matrix $A$ into $p$ blocks of almost equal size, where $p$ is the number of cores.
For cache reuse, the 
blocks assigned to
each core are further partitioned.
Within each block, an orthotrope-style~\cite{PITCH_temporal_overview} temporal blocking is used 
to perform MPK computation locally for the block. This requires to find neighbors of each block that are involved in an MPK computation with power $k$.
However, 
these neighbors end up 
in nonconsecutive spots, resulting in a performance bottleneck.  In~\cite{schreiber_multicore}, MPK kernels 
were studied on modern multicore architectures for banded sparse matrices that arise from stencil discretization.
Following classic stencil blocking approaches, a geometrical 
	blocking method was proposed.
For matrices arising from two-dimensional discretization the method achieved decent speedup.
However, for matrices from three-dimensional discretization it yielded very limited performance gains due to high matrix bandwidth.
Most of the other works~\cite{demmel_MPK_tr,Ichi_GMRES_paper,Ichi_precon} on MPK schemes focused on reducing the MPI communication overhead. 
A recent work~\cite{diamond_MPK} in this direction presents a theoretical study on the benefit of diamond tiling for reducing communication.

\subsection*{Contribution and Outline}

Our work 
bridges the gap
between  temporal blocking of stencil algorithms~\cite{doi:10.1137/070693199,6012879,9355259}, which can be considered 
as an MPK on structured grids,  and recursive 
spatial 
blocking strategies for SpMV~\cite{10.1145/3293883.3295712}. 
In addition we reduce the need to manually set up the blocks. 
We cover full thread-level parallelization and focus on a single multicore processor. Our contributions are as follows:
\begin{itemize}



\item We generalize temporal tiling strategies known from stencil computations 
on structured grids to MPK computations on structured and unstructured sparse matrices using the levels of the graph of the matrix.

\item We present an efficient, multi-threaded implementation of our level-based blocking method for sparse MPK on modern multicore processors. Our solution aims to reduce the main memory traffic and to avoid scalability bottlenecks such as synchronization overhead or load imbalance. 

\item We conduct a detailed performance analysis of our approach as implemented in RACE on various CPU architectures. 

\item For a broad set of sparse matrices we demonstrate full threading functionality and excellent multicore performance achieving speedups of 3$\times$ to 5$\times$ compared to a standard baseline implementation.

\item We validate the performance improvements using the \rlm\ and the phenomenological Execution-Cache-Memory (ECM) model. These models corroborate the optimality of both our level-blocking approach and the baseline implementation to which we compare. 


\item We finally discuss potential applications of MPK to problems from physics and chemistry and demonstrate how such applications can be optimized using the MPK model.
\end{itemize}


The remainder of the paper is structured as follows. 
Section \ref{sec:hardsoft} reviews our experimental  setup, in particular hardware and software characteristics of 
the next generation of scalable processor, namely the  Intel Cascade Lake and Intel Ice Lake, and
 the AMD EPYC architectures, and, additionally, the set of benchmark matrices.
In Section \ref{sec:MPK_kernel} we review the computational workload of matrix-vector 
multiplications for sparse matrices.
Section~\ref{sec:level_based_mpk} is dedicated to the main contribution of the
paper and describes in detail the algorithmic components of 
level-based blocking of MPK.
Section \ref{sec:param_study} includes an assessment of performance parameters within our recursive 
level-based blocking engine (RACE MPK) method.
In Section \ref{sec:perf_eval} we conduct a detailed performance analysis of our cache-aware 
implementation for matrix-power kernels and compare it to a state-of-the-art implementation.
Section \ref{sec:app_cheb_tp} presents the application of the matrix-power-vector multiplication
in Chebyshev time propagation of quantum wave functions and finally Section~\ref{sec:conclusion} concludes the paper.

\section{Hardware and software environment}
\label{sec:hardsoft}

\subsection{Hardware}
\label{sec:testbed}
The measurements in this paper were conducted on a single socket of Intel Cascade Lake (CLX), Intel Ice Lake (ICL), and AMD Epyc Zen2 (ROME), respectively. Key specifications  of the three systems are summarized in Table~\ref{tab:testbed}.

These state-of-the-art processors power more than 50\% of the top 100 ranking supercomputers~\cite{top500_list}.
The Intel CPUs support the \mbox{AVX-512} instruction set, while
the AMD CPU supports only \mbox{AVX-2}. 
Turbo mode was active for all the runs, and
the systems were configured with one ccNUMA domain per socket, i.e., on Intel systems the Sub-NUMA Clustering (SNC) was disabled and on AMD the NPS1 mode was used.

All CPUs have three levels of cache: private, inclusive L1 and L2, and a victim-type L3. 
The L3 cache on the Intel systems is shared by all cores of a socket, while 
on ROME it is shared only within a core complex (CCX), which comprises four cores.
The aggregate L3 cache on ROME is 2.5$\times$ larger than on ICL and 5$\times$ larger than on CLX. This can be observed in the full-socket load-only bandwidth measurements in Fig.~\ref{fig:bandwidth}, where the combined L2 and L3 cache sizes are marked with dashed lines.
This data also shows the L3 and main memory bandwidths of the three CPUs. 
CLX and ICL have a moderate L3 bandwidth of 300\,\GBS\ and 400\,\GBS, respectively, while ROME has a very high L3 bandwidth of more than 2500\,\GBS.
It is worth noting that the transition from L3 to main memory is very sharp on ROME and occurs exactly where the data-set size exceeds the total cache size,
while on the Intel systems the drop is gradual and there is a noticeable cache effect even when the working set exceeds the cache size by 2$\times$ or more, due to its dynamic cache replacement policy~\cite{Understanding}.
The main memory bandwidth ($b_\mathrm{Mem}$) of CLX, ICL and ROME is about 116\,\GBS, 170\,\GBS, and 146\,\GBS, respectively.

\subsection{Software}
For compilation, Intel compilers (see Table~\ref{tab:testbed} for version info) were used on
Ubuntu 18.04.5 (CLX and ROME) and Red Hat Enterprise Linux 8.2 (ICL), respectively,
with compiler flags \texttt{-O3 -xHOST}.
All floating-point computations were done in double precision, while integers were 32~\bits\ wide.
Threads were bound to cores in a closed (fill-type pinning) manner.
To reduce fluctuations, each kernel was executed multiple times 
such that the overall runtime is greater than one second.
The average performance of these runs was then reported.
As the variation among multiple measurements was less than 5\%, we do not show error bars.

\begin{figure*}[tbp]
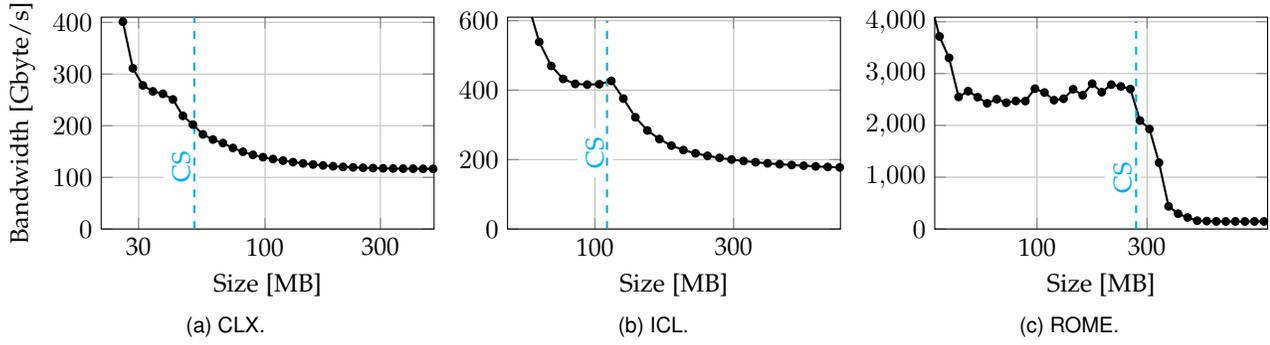

	\subfloat[CLX.]{%
	\input{plots/tikz/bandwidth/casclakesp2/load_avx512.tex}%

		\label{fig:bandwidth:clx}}
	\subfloat[ICL.]{%
	\input{plots/tikz/bandwidth/horeka/load_avx512.tex}%

	\label{fig:bandwidth:icl}}
	\subfloat[ROME.]{%
	\input{plots/tikz/bandwidth/tg094/load_avx.tex}%

		\label{fig:bandwidth:rome}}
	\caption{Single socket LLC and memory bandwidth (load-only) of the three architectures under consideration. The dashed line represents the total available cache size (CS). 
	Note the different scaling on the $y$-axis.\label{fig:bandwidth}}
\end{figure*}

For pinning, bandwidth benchmarks (see Fig.~\ref{fig:bandwidth}),
and for counting hardware events
we use the \likwidpin, \likwidbench, and \likwidperfctr\ tools from the \likwid\ tool suite version 5.1.
\begin{table}[!tb]
	\centering\footnotesize
	\caption{\label{tab:testbed}Key specification of test bed machines.}
	\resizebox{\linewidth}{!}{%
	\begin{tabular}{l c c c}
		\toprule
		Architecture  & CLX  & ICL  & ROME  \\
		\midrule
		Chip Model                      & Xeon Gold 6248                & Xeon Platinum 8368 & AMD EPYC 7662 \\
		Microarchitecture  & Cascade Lake  & Sunny Cove  & Zen-2 \\
		Release year & 2019 & 2021 & 2020\\
		Cores per socket  & 20     &    38 & 64   \\
		Max. SIMD width       & 512\,\bits & 512\,\bits & 256\,\bits\\
		L1D cache capacity   & 20$\times$32\,\KiB   & 38$\times$48\,\KiB   & 64$\times$32\,\KiB   \\
		L2 cache capacity    & 20$\times$1\,\MiB    & 38$\times$1.25\,\MiB  & 64$\times$512\,\KiB     \\
		L3 cache capacity   & 27.5\,\MiB &  57\,\MiB    & 16$\times$16\,\MiB \\
		Memory Configuration &  6 ch. DDR4-2933  & 8 ch. DDR4-3200 &  8 ch. DDR4-3200   \\
		Mem. Bandwidth ($b_\mathrm{Mem}$)      & 116\,GB/s  & 170\,GB/s & 146\,GB/s   \\
		Operating system   & Ubuntu 18.04.5  & RHEL 8.2 & Ubuntu 18.04.5    \\ 
		Compiler	& Intel 19.0 update 5 & Intel 19.1 update 3 & Intel 19.0 update 5 \\
		\bottomrule
	\end{tabular}
	} 
\end{table}


\subsection{Benchmark matrices}
\begin{table}[tbp]
	\caption{Details of the benchmark matrices. \NR\ is the number of
		rows, \NNZ\ is the number of nonzeros, and \NNZR\
		is the average number of nonzeros per row.
		\label{tab:matrices}}
	\begin{center}
		\resizebox{\linewidth}{!}{%
			\begin{tabular}{|l|l|S[table-format=7.0, table-space-text-pre=(, table-space-text-post=)]|S[table-format=8.0, table-space-text-pre=(, table-space-text-post=)]|S[round-mode=places,round-precision=2]|}
\toprule
{Index} & {Matrix name} &  {\NR} & {\NNZ} & {\NNZR} \\
\midrule
{1} & {cfd2} & 123440 & 3087898 & 25.01537589 \\
{2} & {parabolic\_fem} & 525825 & 3674625 & 6.988304094 \\
{3} & {xenon2} & 157464 & 3866688 & 24.5560128 \\
{4} & {cant} & 62451 & 4007383 &  64.16843605 \\
{5} & {offshore} & 259789 & 4242673 & 16.3312265 \\
{6} & {Hamrle3} & 1447360 & 5514242 & 3.809862094 \\
{7} & {bmw7st\_1} & 141347 & 7339667 & 51.92658493 \\
{8} & {G3\_circuit} & 1585478 & 7660826 & 4.831871524 \\
{9} & {shipsec1} & 140874 & 7813404 & 55.46377614 \\
{10} & {ship\_003} & 121728 & 8086034 & 66.4270669 \\
{11} & {thermal2} & 1228045 & 8580313 & 6.986969533 \\
{12} & {gearbox} & 153746 & 9080404 &  59.06107476 \\
{13} & {crankseg\_1} & 52804 & 10614210 &  201.0114764 \\
{14} & {pwtk} & 217918 & 11634424 & 53.38899953 \\
{15} & {rajat31} & 4690002 & 20316253 & 4.33182182 \\
{16} & {gsm\_106857} & 589446 & 21758924 & 36.91419401 \\
{17} & {F1} & 343791 & 26837113 & 78.06229075 \\
{18} & {cage14} & 1505785 & 27130349 & 18.01741218 \\
{19} & {Fault\_639} & 638802 & 28614564 & 44.79410522 \\
{20} & {inline\_1} & 503712 & 36816342 & 73.09006337 \\
{21} & {RM07R} & 381689 & 37464962 & 98.15572888 \\
{22} & {Emilia\_923} & 923136 & 41005206 & 44.41946365 \\
{23} & {ldoor} & 952203 & 46522475 & 48.85772782 \\
{24} & {af\_shell10} & 1508065 & 52672325 & 34.927092 \\
{25} & {HPCG-128-128-128} & 2097152 & 55742968 & 26.58031845 \\
{26} & {Hook\_1498} & 1498023 & 60917445 & 40.66522677 \\
{27} & {Geo\_1438} & 1437960 & 63156690 & 43.92103396 \\
{28} & {Serena} & 1391349 & 64531701 & 46.38067156 \\
{29} & {bone010} & 986703 & 71666325 & 72.63211422 \\
{30} & {audikw\_1} & 943695 & 77651847 & 82.28489819 \\
{31} & {channel-500x100x100-b050} & 4802000 & 85362744 & 17.77649813 \\
{32} & {dielFilterV3real} & 1102824 & 89306020 & 80.97939472 \\
{33} & {nlpkkt120} & 3542400 & 96845792 & 27.33903342 \\
{34} & {ML\_Geer} & 1504002 & 110879972 & 73.7232876 \\
{35} & {Flan\_1565} & 1564794 & 117406044 & 75.02971254 \\
{36} & {stokes} & 11449533 & 349321980 & 30.50971424 \\
\bottomrule
\end{tabular}

		}
	\end{center}
\end{table}
Table~\ref{tab:matrices} shows the sparse matrices used for the benchmarks and some of their properties:
\NR\ is the total number of rows, \NNZ\ is the total number of nonzero entries, and \NNZR\ is the average number of nonzero entries per row (i.e.,\NNZ/\NR).
The matrices are ordered (top to bottom) according to increasing \NNZ, and all are square since this is a requirement for the matrix power kernel (MPK).
All matrices except one were taken from the SuiteSparse Matrix Collection~\cite{UOF}.
\texttt{HPCG-128-128-128} is the matrix found in the  HPCG benchmark~\cite{HPCG}, with a problem size of $128^3$.

\section{Matrix Power Kernel}
\label{sec:MPK_kernel}
The basic algorithmic workload addressed in this paper is the computation of powers of a sparse matrix applied to a dense vector. The matrix power kernel (MPK) is defined as follows: For a given square, sparse matrix $A$ and a dense vector $x$ calculate all matrix powers $A^px$ up to a maximum $p_m$ ($p=1,\ldots,p_m$) and store all $p_m$ resulting vectors ($y_p=A^p x$) for subsequent calculations. We further define $y_0 := x$.  

\subsection{Baseline MPK implementation}
\label{sec:baseline_MPK}
The standard approach to implement the MPK is to perform a sequence of $p_m$ \spmv\ operations, i.e., $y_i=A y_{i-1}$ with $i=1,\ldots,p_m$, using standard  \spmv\  implementations or library calls. We refer to this strategy as \emph{baseline} MPK.
\begin{figure}[tbp]
	\centering
	\resizebox{\linewidth}{!}{%
		\begin{minipage}{1.34\linewidth}
			\begin{algorithmic}[1]	
				\STATE{$double:: val[$\NNZ$]$} \textcolor{gray}{//store values of nonzeros in $A$}
				\STATE{$int:: col[$\NNZ$], rowPtr[$\NR$+1]$} \textcolor{gray}{//column index and row pointer of $A$}
				\STATE{$double:: y[$\NR, 0:$p_m]$} \textcolor{gray}{//to store input and output vectors}
				\STATE \textcolor{gray}{//Perform $p_m$ SpMVs}
				\FOR{$p=1:p_m$}
					\STATE $y[:,p]$=SpMV($y[:,p-1]$, 0 , \NR-1) 
				\ENDFOR
				\STATE \textcolor{gray}{//Perform SpMV between $A$ and $in$ vector and store result in $out$.\\
				// $row\_s$ and $row\_e$ arguments are used to specify the start and end row\\
				//to which SpMV is applied.}
				\STATE \textbf{function} SpMV($double:: in\_rhs[$\NR$]$, $int:: row\_s$, $int:: row\_e$)
				\begin{ALC@g}
					\STATE{$double::out\_lhs[$\NR$]$} 
					\STATE \textcolor{gray}{//Loop over rows}
					\STATE{\textcolor{darkgray}{\#pragma omp parallel for schedule(static)}}
					\FOR{$row=row\_s:row\_e$}
					\STATE{$double::tmp=0$}
					\STATE \textcolor{gray}{//Loop over nonzeros in row}
					\FOR{$idx=rowPtr[row]:(rowPtr[row+1]-1)$}
					\STATE{$tmp \mathrel{+}= val[idx]*in\_rhs[col[idx]]$} 
					\ENDFOR
					\STATE{$out\_lhs[row] = tmp$}
					\ENDFOR
					\STATE \textbf{return} $out\_lhs$
				\end{ALC@g}
				\STATE \textbf{end function}				

			\end{algorithmic}
		\end{minipage}
	}
	\caption{CRS-based MPK computing $A^{p_m}x$. 
		The arrays $val$, $col$, and $rowPtr$ hold the CRS data structure of $A$.
		The input and output vectors are stored in the $y$ matrix.
	}
	\label{fig:SpMV_alg}
\end{figure}

\begin{figure*}[tbp]
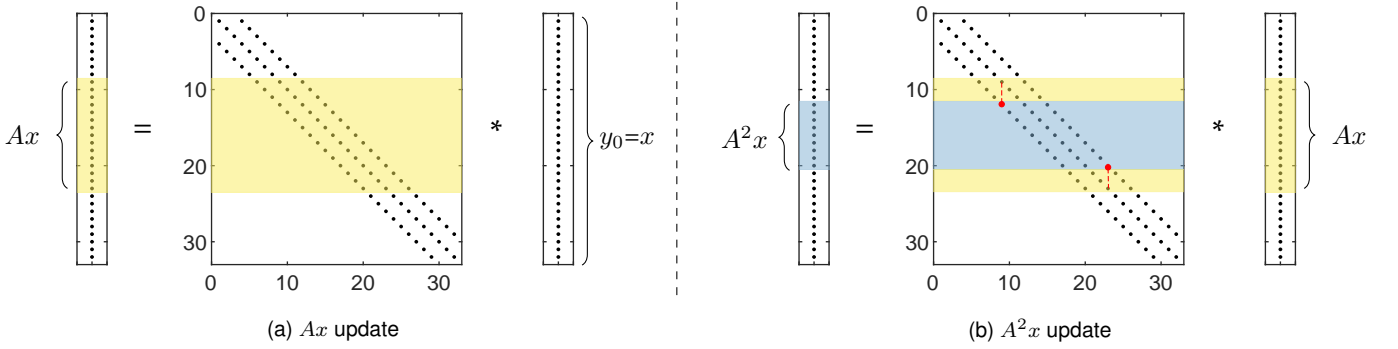

	\centering	
	\subfloat[$Ax$ update]{%
	\input{plots/tikz/simple_diag_example/simple_diag_p1.tex}%

		\label{fig:simple_diag:p1}}
	\subfloat[$A^2x$ update]{%
	\input{plots/tikz/simple_diag_example/simple_diag_p2.tex}%

		\label{fig:simple_diag:p2}}
	\caption{Blocking successive matrix applications for a simple banded sparse matrix: 
		(a) The RHS vector is the input vector $x$. Yellow elements of the LHS vector are updated to $Ax$. (b) The next update is performed on the blue block of $A$ to compute $A^2x$ on the blue elements of the LHS vector.	
		Yellow matrix elements can be reused when computing $A^2x$  on blue blocks. 
		\label{fig:simple_diag}}
\end{figure*}
Figure~\ref{fig:SpMV_alg} shows a high-level representation of our baseline MPK together with an \spmv\ implementation that is known to provide good performance on CPUs for a wide variety of sparse matrix structures. 
The sparse matrix $A$ is stored in the well-known CRS format,
using the three arrays $rowPtr$, $val$, and $col$, which hold the 
row pointer information, values, and column indices of nonzero entries, respectively (see~\cite{CRS_saad} for details).
This information is passed (as global data) to the \spmv\ function along with the function parameters (line 9) representing the right-hand side (RHS) vector and the range of row indices for which the SpMV is to be computed.\footnote{For the baseline implementation, the entire row range is specified.}
The function then performs the \spmv\ operation (lines 12--20) and returns the resulting left-hand side (LHS) vector.
Note that most \spmv\ implementations in libraries are unsuitable for the optimized MPK discussed later as they do not support \spmv\ on a subset of rows. 
Therefore  we use our own version of \spmv, which serves as the main kernel for both the baseline and the optimized version.
We have ensured that our \spmv\ performs at least as good as Intel MKL with the standard CRS format.

The baseline MPK stores the $p_{m}+1$ vectors $\{y_p\}$ in the matrix $y[:,0:p_m]$ (column-major order) 
and performs $p_m$ back-to-back calls to the \spmv\  
function (see lines 5--7 of Fig.~\ref{fig:SpMV_alg}).  
If the caches are too small to hold the entire matrix, it must be read $p_m$ times from main memory.
Consequently, the optimum (minimum) main memory balance for the CRS-based baseline MPK is $B_C=6\,\BF$ \cite{Gropp99towardsrealistic,RACE}, which is equivalent to $12\,\bytes$ of memory traffic per nonzero matrix entry. 
The baseline MPK thus reflects the strongly memory-bound performance characteristic of the underlying \spmv\ operation.

In order to evaluate the quality of optimized MPK implementations, we will measure the actual code balance $B_{C,m}$
and compare it with the theoretical baseline minimum ($6\,\BF$) discussed above.
The $B_{C,m}$ is obtained by measuring the actual data traffic (using \likwidperfctr) and dividing it by the minimum amount of floating-point operations to be performed, i.e., $2\times N_{nz}\times p_{max}$. Where appropriate, measured code balance from within the cache hierarchy will also be reported.

\subsection{Blocking strategy for the MPK implementation}
\label{sec:blocking_basic_idea}
As the same sparse matrix is repeatedly applied, there is substantial performance optimization  potential via data transfer reduction by reusing matrix entries from the cache for the successive computation of multiple powers.
The basic idea is to compute  the \spmv\ partially for a block of $A$ that fits into cache and reuse these matrix entries for the next \spmv, i.e., calculate another power on a smaller subset of the data.
This approach is equivalent to temporal blocking for iterative stencil update schemes, where multiple updates on the same stencil data are computed in cache. 
Here the spatial stencil structure determines the dependencies between successive updates and geometric schemes for handling the spatial-temporal dependencies such as trapezoidal~\cite{Frigo} or diamond blocking~\cite{diamond_blocking} are well established. 
To demonstrate the equivalent challenge in MPK, we show in Fig.~\ref{fig:simple_diag} a simple banded sparse matrix, which arises from a discretization of a toy stencil in one spatial dimension.
In the first step (Fig.~\ref{fig:simple_diag:p1}), an \spmv\ operation is performed applying a block of the matrix (yellow rows), which fits into cache, to the input (RHS) vector $x$ to calculate a part of $A x$ (yellow elements of LHS vector).
In the next step (Fig.~\ref{fig:simple_diag:p2}), the updated vector elements serve as input and are used to calculate $A^2 x$ (blue elements of the LHS vector) by applying \spmv\ with a subset of the matrix block (blue rows). 
To fulfill the dependencies between these successive \spmv\ steps, the column indices of the subset of the matrix block (blue rows in Fig.~\ref{fig:simple_diag:p2}) must be in the range (indicated with red line) of the row indices of the original matrix block (yellow rows). 
It is obvious that the overhead of this approach, which is quantified by the ratio of yellow to blue rows in Fig.~\ref{fig:simple_diag:p2}, increases with the bandwidth of the matrix (i.e., with longer-range stencils). 

The outlined MPK blocking approach can be generalized for sparse matrices with irregular structures. 
We define ${\cal I}$ to be a set of row indices of the matrix $A$. The corresponding set ${\cal C}({\cal I})$ contains the column indices of all nonzero entries in the rows of $\cal I$, i.e., if $i \in \cal I$ then $j\in {\cal C}({\cal I}) \iff A_{i,j} \ne 0$.
Based on this notation, the \spmv\ operation ($y=Ax$) for a given row index $i \in {\cal I}$ can be written as:
\begin{equation}
\label{eq:SpMV_index}
	y_i = \sum\limits_{j \in {\cal C}(i)}^{} A_{i,j} x_j
\end{equation}
If we apply the \spmv\ for all rows in $\cal I$ to a RHS $y_{p-1}$, then all corresponding row entries of the LHS vector are updated to power $p$. We can then apply to this vector another \spmv\ on a set of rows $\cal K$ for which ${\cal C}({\cal K}) \subseteq \cal I$.

The choice of the set of row indices $\cal I$ for a given sparse matrix $A$ is decisive to the performance of such a method:
(i) The matrix elements associated with $\cal I$ and ${\cal C}({\cal I})$ have to fit into cache and (ii) should be stored 
to enable high spatial and temporal locality without
indirect access.
Furthermore, (iii) the bandwidth of the matrix involved in the MPK should be as small as possible, i.e., the indices of ${\cal C}({\cal I})$ have to be close to the set $\cal I$.
A potential approach to address these challenges is to consider the \spmv\ operation as a graph traversal problem as done in the RACE coloring scheme~\cite{RACE}. 
Here, breadth-first search (BFS)~\cite{BFS} is applied to the graph underlying $A$; the BFS levels of $A$ are stored consecutively. 
These levels allow us to identify appropriate parts of the matrix ($\cal I$ and $\cal K$)  for blocking 
and how to traverse the full matrix (graph) systematically to update vector elements corresponding to all matrix powers while maintaining locality in accessing matrix and vector data.    
As an added benefit, the BFS reordering of the matrix reduces its bandwidth.

\section{Level-blocked MPK}
\label{sec:level_based_mpk}
The RACE coloring scheme has been developed to generate hardware efficient distance-$k$ colorings of graphs~\cite{RACE}. 
It has been successfully applied to the shared-memory parallelization of symmetric \spmv\ providing unprecedented performance levels.
Further it has been shown that the level-based approach allows to control dependencies in \spmv\ operations and at the same time provides flexibility to ensure data locality and to adjust to the degree of parallelism required by modern multicore processors. 
The level-based blocking strategy introduced in the following adapts these RACE properties to the MPK. 
We thus first recapitulate the basic terminology and the level-based approach of RACE. 
Next we demonstrate how it is basically applied to the MPK problem and then show how data locality and efficient shared-memory parallelization can be achieved. 

In this section, we restrict ourselves to symmetric matrices, i.e., undirected graphs. 
However, the proposed MPK blocking method is also applicable to non-symmetric square matrices. 
The following definitions from graph theory are used throughout the paper:\\[2mm]
\textbf{Graph:} $G = ({\cal V}, {\cal E})$ represents a graph, with ${\cal V}(G)$ denoting a set of vertices and ${\cal E}(G)$ denoting its edges. 
	For sparse matrices,  ${\cal V}(G)$ consists of all row indices of the matrix and ${\cal E}(G)$ consists
	of edges between two vertices corresponding to the row ($u$) and the column indices ($v$) of the nonzero entries, i.e., $\{u,v\} \in {\cal E}(G) \iff A_{u,v} \ne 0$. 
	\\	
\textbf{Neighborhood:} The neighborhood of a vertex $u$ is the set of vertices ${\cal N}(u)=\{v \in {\cal V}(G): \{u,v\} \in {\cal E}(G)\}$.\\
\textbf{Subgraph:} A subgraph $H$ of $G$ specifically refers to the subgraph induced by
	vertices ${\cal V}' \subseteq {\cal V} (G)$ and is defined as the graph 
	$H= ({\cal V}',\{\{u,v\}\in {\cal E}(G) \land  u,v \in {\cal V}'\})$.\\[1mm]

In the graph terminology, an \spmv\ operation ($y=Ax$) can be formulated as follows: 
If  $G = ({\cal V}, {\cal E})$ is the graph representation of the sparse matrix $A$ 
 then for every vertex $u \in {\cal V}(G)$ calculate
 \begin{equation}
\label{eq:SpMV_graph}
	y_u = \sum\limits_{v \in {\cal N}(u)}^{} A_{u,v} x_v\eos
\end{equation}
Comparing ~\eqref{eq:SpMV_graph} with~\eqref{eq:SpMV_index}, we can observe the equivalence between index-based (row index $i$ and its related column indices ${\cal C}(i)$) and graph-based (vertex $u$ and its neighborhood ${\cal N}(u)$) notations.

%
%
\begin{figure*}[tbp]
	\subfloat[Graph.]{%
	\input{plots/tikz/2d-7pt/graph_2d_7pt.tex}%

		\label{fig:stencil_graph}}
	\subfloat[Sparsity pattern of matrix.]{\includegraphics{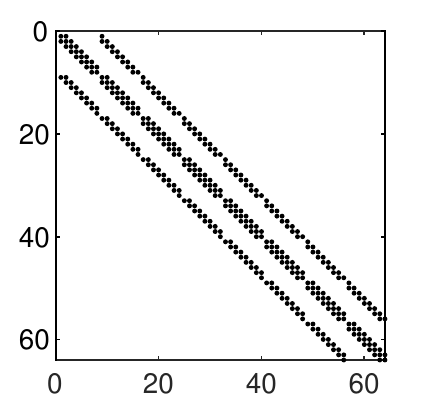}\label{fig:stencil_matrix}}
	\subfloat[Permuted graph.]{%
	\input{plots/tikz/2d-7pt/graph_2d_7pt_permuted_w_levels.tex}%

		\label{fig:stencil_graph_permuted}}
	\subfloat[Permuted matrix.]{\includegraphics{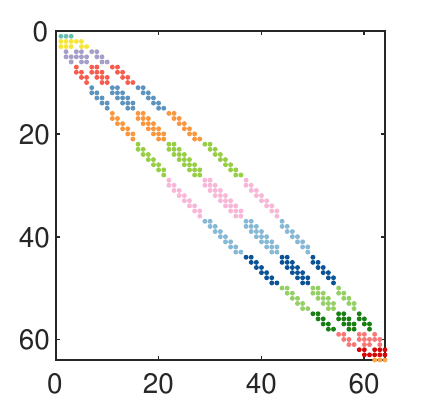}\label{fig:stencil_matrix_permuted}}	
	\subfloat[]{\includegraphics{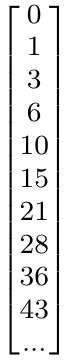}\label{fig:level_ptr}}
	\caption{
		Graph~(a) and sparsity pattern~(b) of the matrix associated with a 2d-7pt stencil on an  8$\times$8 grid. 
		In~(a), the associated stencil is highlighted in red for an arbitrary vertex (54).
		(c) shows the permuted graph and (d) the sparsity pattern of the matrix after applying BFS reordering.
	The vertices (rows) of the graph (matrix) that belong to a level are 
represented with the same color.
The \texttt{level\_ptr} associated with the permuted graph/matrix is shown in~(e)}
\end{figure*}
%
To illustrate our method, a simple graph generated by applying a two-dimensional seven-point (2d-7pt) stencil to a square grid of size 8$\times$8 will serve as an example.
Figure~\ref{fig:stencil_graph} shows the graph with each vertex
numbered in lexicographic ordering. 
The associated stencil at a single grid point (vertex 54 and its neighborhood) is highlighted.
The sparsity pattern of the corresponding matrix is shown in
Fig.~\ref{fig:stencil_matrix}.

\subsection{Levels}
\label{sec:levels}
The level formation in RACE is based on a BFS which assigns each vertex (row) of the graph (matrix) to a level.
First, a \emph{root vertex} $v_\mathrm{root}$ is chosen and assigned to the first level, $L(0)$.
The rest of the levels, $L(i)$ $\forall$ $i>0$, are defined to contain 
vertices that are in the combined neighborhood of the vertices in the previous level $L(i-1)$ but have no level numbers assigned yet, i.e.,
\begin{equation}\label{eq:level}
L(i) = 
\begin{cases}
v_\mathrm{root} & \text{ if } i = 0, \\
\big\{ u:  u \in {\cal N}(L(i-1))\, \wedge \\
\textcolor{white}{u:\big\{}  u \not\in \{L(0), \dots, L(i-1)\}\big\}  & \text{ else}.\\
\end{cases}   
\end{equation}
Figure~\ref{fig:stencil_graph_permuted} shows the 15 levels (indicated by different colors) generated by this procedure for the stencil graph if $v_\mathrm{root}=0$ is chosen.
After level formation, the vertices are renumbered (compare vertex indices in Fig.~\ref{fig:stencil_graph} and Fig.~\ref{fig:stencil_graph_permuted}) such that those in the same level are numbered consecutively and the vertices in level $L(i-1)$ appear before those in $L(i)$.
This permutation\footnote{Note that a symmetric permutation is employed on the matrix, i.e., both rows and columns are permuted} 
increases data locality between neighboring vertices and results in a lens-shaped matrix with typically reduced bandwidth (see Fig.~\ref{fig:stencil_matrix_permuted}).
Since this improves the data locality of sparse matrix computations, such permutations are widely employed as  preprocessing steps for \spmv-based algorithms~\cite{RCM_SpMP}. 

As a consequence of the definition of levels, the neighborhood of all vertices in a given level $L(i)$ is clearly confined to the vertices within the previous, current, and next levels, i.e.:
%
\begin{equation}\label{eq:level-neighb}
{\cal N}(L(i)) \in \{L(i-1) \cup L(i) \cup L(i+1)\}, \text{ for } i>0\eos
\end{equation}
This property is crucial for the design of our level-based MPK blocking scheme as it defines the dependency between the computation of {\spmv}s for different levels at different matrix powers:
To advance all vertices of $L(i)$ to $A^px$, the calculation of $A^{p-1}x$ has to be completed on the levels $L(i-1)$, $L(i)$, and $L(i+1)$.




\subsection{Level-based blocking of \MPK}
\label{sec:level_blocking_idea}
In Sec.~\ref{sec:MPK_kernel}, we discussed the baseline MPK and the 
potential matrix data reuse by blocking across the \spmv\ operations involved in the MPK.  
Further it has been shown that the graph formulation of the \spmv~\eqref{eq:SpMV_graph} together with the neighborhood relation~\eqref{eq:level-neighb} of the levels~\eqref{eq:level} provide a natural framework for the structured computation of the MPK. 
This includes the dependency between a level and its neighborhood; e.g., in Fig.~\ref{fig:stencil_graph_permuted} one can calculate the next matrix power for level $L(6)$ (with vertices $21,\ldots,27$) only after the computation of the previous matrix power is complete on levels $L(5)$, $L(6)$, and $L(7)$ (containing vertices $15,\ldots,35$).

\begin{figure}[tbp]
	\centering
	%
	\input{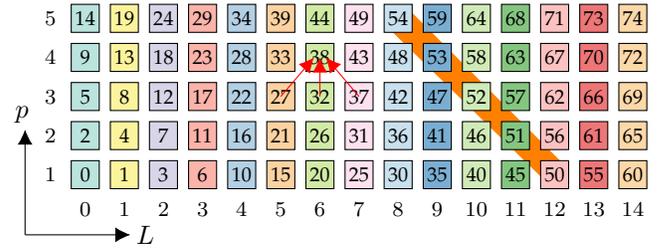}%

	\caption{$Lp$ diagram with 15 levels ($L(0),\ldots,L(14)$) and a maximum power stage of $p_{max}=5$. Level colors are the same as in Fig.~\ref{fig:stencil_graph_permuted}.
	Each node in the $Lp$ diagram is numbered according to the execution order.
	For $p=4$ and level $L(6)$, the explicit dependencies with levels at $p=3$ are indicated with red arrows. 
	The nodes highlighted in orange fulfill $i+p=13$ (``diagonal'').
	\label{fig:lp_basic}}
\end{figure}

We next introduce the \emph{$Lp$ diagram} to visualize the dependencies between levels in MPK calculations.
In the $Lp$ diagram, the indices of the levels $L(i)$ are on the $x$-axis and the matrix power stages ($1\leq p \leq p_{max}$) are on the $y$-axis.
Hence, each node $(i,p)$ in the diagram represents an \spmv\ on the vertices in level $i$ to compute part of the power $p$\@.
Figure~\ref{fig:lp_basic} shows the $Lp$ diagram for 15 levels and $p_{max}=5$.
To satisfy the dependencies in the level-based MPK blocking scheme, the nodes $(i-1,p-1)$, $(i,p-1)$, and $(i+1,p-1)$ need to be computed  before \spmv\ can be applied to compute the node $(i,p)$.
The red arrows in Fig.~\ref{fig:lp_basic} denote the dependency for the computation of $L(6)$ at $p=4$, i.e., for the node $(6,4)$.
The order of traversal in the $Lp$ diagram is as follows:
\begin{itemize}
  \item Each \emph{diagonal}, defined by $i+p=\text{const}$, is traversed from bottom to top (starting at $p=1$).
  \item Diagonals are traversed from left to right, i.e., starting with $p=1$ for $L(0)$. 
\end{itemize}
This execution order, which is independent of the actual graph structure,
ensures that the levels $L(i-1)$, $L(i)$, and $L(i+1)$ are updated to power stage $p-1$ before level $L(i)$ is advanced to power stage $p$.
In Fig.~\ref{fig:lp_basic}, the order of all execution steps of this scheme is shown via the node numbers
in the $Lp$ diagram
 with $p_m=5$.

Visualizations similar to Fig.~\ref{fig:lp_basic} are often shown for one-dimensional (1D) radius-one stencils, where the $x$-axis represents the grid points and the $y$-axis shows iterations or time steps%
~\cite{Frigo,diamond_blocking,PITCH_temporal_overview}. 
As we have shown above, our level-based MPK algorithm shows the same dependencies, with levels substituting grid points on the $x$-axis.
This opens up a host of options, since we could draw from the large variety of temporal blocking optimizations developed for 1D stencils. 
Our approach is analogous to parallelogram-style temporal blocking; see~\cite{PITCH_temporal_overview} for a classification.

The reuse distance of a given level is a central quantity to characterize the cache locality of the level-blocked (LB) MPK\@. 
Within the $Lp$ diagram, this quantity can be determined by the number of execution steps between two computations on the same level, i.e., one step in vertical direction. 
As the scheme traverses the $Lp$ space in consecutive diagonals, 
a level computed at power $p$ will be reused after 
$\tilde{d} + 1$ execution stages for the computation of the next power $p+1$, where $\tilde{d}$ is the number of execution steps in the current diagonal. 
After the wind-up and before the wind-down phases at the left and right ends of the $Lp$ diagram, we have $\tilde{d}=p_m$; hence, levels are reused after $p_m+1$ execution steps.
This can be observed from Fig.~\ref{fig:lp_basic}
if we concentrate on a single level, e.g., the vertices of $L(10)$ 
used in the 40th execution step to compute $p=1$ are reused in the 46th step to compute $p=2$.
As the number of levels is typically much larger than the maximum power stage, we can assume a maximum reuse distance of $p_m+1$ execution stages.
This means if all the matrix entries associated with the $p_m+1$ successive levels touched 
between two computations of a given $L(i)$ can be held in a cache, all accesses to this $L(i)$ can be served from the cache with the exception of the first one ($p=1$), which requires main memory access. 
Assuming that cache accesses are much faster than memory accesses, the performance of the LB MPK implementation can improve by a factor of at most $p_m$ as compared to the baseline MPK.
\begin{figure}[tbp]
	\centering
	\resizebox{\linewidth}{!}{%
		\begin{minipage}{1.34\linewidth}
			\begin{algorithmic}[1]
				\STATE{\textcolor{gray}{//traverse diagonals of $Lp$ in ascending order $d=i+p$}}
				\FOR{$d = 1:L_m+p_m-1$} 
				\STATE $p_\mathrm{start} = \max(1,$ $d-(L_m-1))$
				\STATE $p_\mathrm{end} = \min(d,$ $p_m)$
				\STATE{\textcolor{gray}{//traverse diagonal $d=i+p=const$ in ascending order of p}}
				\FOR{$p = p_\mathrm{start}:p_\mathrm{end}$}
				\STATE $i = (d-p)$ \textcolor{gray}{//$i+p=d$ diagonal}
				\STATE $y[p]$ = SpMV($y[p-1]$, $level\_ptr[i]$, $level\_ptr[i+1]-1$) 
				\ENDFOR
				\ENDFOR
			\end{algorithmic}
		\end{minipage}
	}
	\caption{Basic implementation of the level-blocked (LB) MPK algorithm.
		$L_m$ is the total number of levels and $p_m$ is the maximum matrix power. The \spmv\ function implementation from Fig.\ref{fig:SpMV_alg} is used.
	}
	\label{fig:MPK_basic_alg}
\end{figure}

\begin{figure*}[tbp]
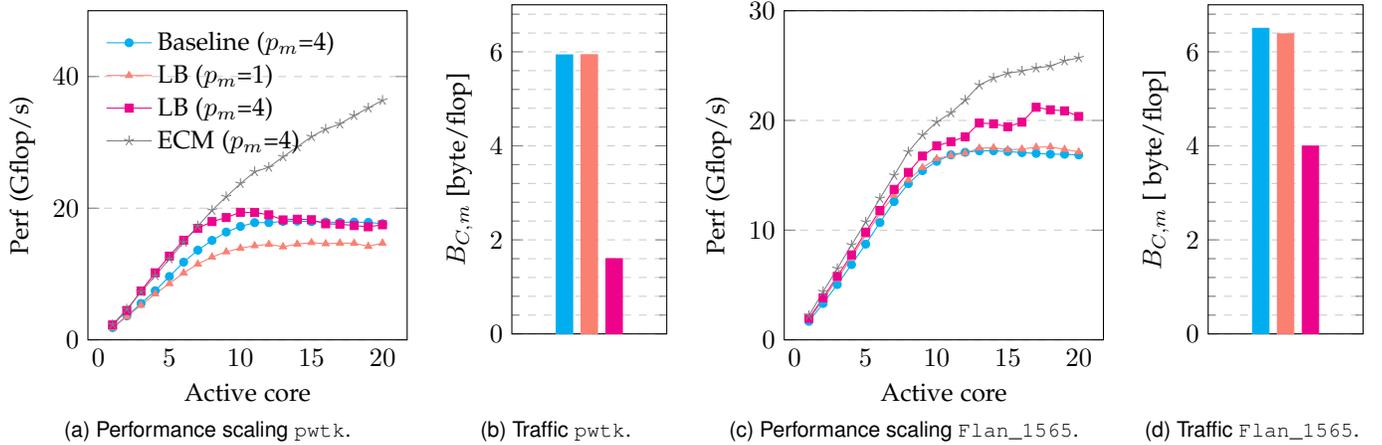

	\centering
	\subfloat[Performance scaling \texttt{pwtk}.]{%
	\input{plots/tikz/mtxPower_results/motivation/plot_barrier_cs1.tex}%
\label{fig:motivation:initial:scaling}}
	\hfill
	\subfloat[Traffic \texttt{pwtk}.]{%
	\input{plots/tikz/mtxPower_results/motivation/plot_barrier_cs1_traffic.tex}%
\label{fig:motivation:initial:traffic}}
	\hfill
	\subfloat[Performance scaling \texttt{Flan\_1565}.]{%
	\input{plots/tikz/mtxPower_results/motivation/plot_barrier_cs1_Flan_1565.tex}%
\label{fig:motivation:initial:scaling:Flan}}
	\hfill
	\subfloat[Traffic \texttt{Flan\_1565}.]{%
	\input{plots/tikz/mtxPower_results/motivation/plot_barrier_cs1_traffic_Flan_1565.tex}%
\label{fig:motivation:initial:traffic:Flan}}
	\caption{Scaling performance and main memory traffic of our LB MPK implementation for $Ax$ ($p_m=1$) and $A^4x$ ($p_m=4$) in comparison to the baseline MPK on one socket of CLX for the \texttt{pwtk} and \texttt{Flan\_1565} matrices.
	  The stars show the phenomenological ECM performance model~\cite{Malas2017} (in gray) for the $p_m=4$ case. 
	  The model assumes that the computation of first power of a level ($p=1$) does not overlap for subsequent powers ($p>1$).
	}\label{fig:motivation:initial}
\end{figure*}


\subsection*{Implementation}
Two basic implementation decisions for our LB MPK are guided by RACE. First, the complete algorithm operates on the permuted graph. Second, only two lean data structures are required to store the information on the permutation and the levels:  
The permutation vector (\NR\ entries) is required to recover the original ordering. 
The storage location of the first vertex (row) of each level are stored in the  \texttt{level\_ptr} array (one entry per level). 
Figure~\ref{fig:level_ptr} shows the \texttt{level\_ptr} of our stencil example matrix (see Fig.~\ref{fig:stencil_matrix_permuted}).

A straightforward implementation of our LB MPK is presented in Figure~\ref{fig:MPK_basic_alg}.
The algorithm first iterates over all diagonals of the $Lp$ diagram in ascending order (line 2).
Within a diagonal $d=i+p=const$, the computations are processed in increasing order of power $p$ (line 6). 
Note that to account for the wind-up and the wind-down phase of the parallelogram,
the starting and ending power stages are adjusted in lines 3 and 4 of the algorithm.
Depending on the power $p$ and the diagonal counter $d$, the actual level index $i$ to use in the current iteration
is calculated in line 7.
Finally, in line 8 the vector ($y[p-1]$) containing the required information at power level $p-1$ and the indices of the first and last row of $L(i)$ are passed to the \spmv\ function (shown in Fig.~\ref{fig:SpMV_alg}) to compute $A^{p}x$ on level $L(i)$. 
Note that OpenMP parallelization is done within the \spmv\ function using static scheduling (line 12 in Fig.~\ref{fig:SpMV_alg}). 
As there is an implicit barrier after the parallel workshare construct, all threads finish the computations on a given execution stage before proceeding to the next one.
In order to reduce the start-up overhead at the parallel region encountered in each \spmv\ call, the parallel region is opened outside
 the \spmv\ routine in our implementation.


Note that the storage of each level is consecutive and the levels are stored in ascending index order.
Therefore, the proposed method neither has irregular accesses to matrix entries nor 
does it have to store extra copies of matrix elements and perform redundant computations,
which were required in previous work~\cite{mohiyudeen}.
Moreover,the parallelization within the levels avoids load imbalance and redundant thread-local copies, which may add significant overhead for irregular matrices and high thread (or core) counts. 

\subsection*{Performance analysis of naive version}
The naive implementation of the LB MPK already results in a decent performance improvement for some of the matrices presented in Table~\ref{tab:matrices}.
However, it often falls short of the predicted maximum $p_m$-fold speedup. 
For example, with $p_m=4$ on one socket of CLX, 50\% of the matrices in the table showed speedup of less than 10\% and almost 10 matrices had a performance degradation compared to the baseline MPK\@.
We choose two representative matrices, \texttt{pwtk} and \texttt{Flan\_1565}, which are exemplary for the major performance shortcomings of the basic LB MPK and we will identify those in the following.

Figure~\ref{fig:motivation:initial} shows the multithreaded performance and main memory code balance of the LB MPK  (Fig.~\ref{fig:MPK_basic_alg}) with $p_m=1$ and $p_m=4$ along with the baseline MPK (Fig.~\ref{fig:SpMV_alg}) with $p_m=4$ on one socket of CLX (20 cores) for both matrices.
One may expect that LB MPK with $p_m=1$ and the baseline MPK  should deliver the same performance, independent of $p_m$. 
They both perform the memory-bound \spmv\ operations successively but with different execution order within each \spmv\ function, and their minimum code balance from main memory is $B_C=6\,\BF$ (see Sec.~\ref{sec:baseline_MPK}).
Hence, a data traffic (i.e. $B_C$) reduction and performance speedup of at most 4$\times$ may be achieved when using LB MPK for $p_m=4$. 

For \texttt{pwtk}, the typical memory bandwidth saturation pattern is observed for LB MPK ($p_m=1$, triangles) and baseline MPK (circles) in Fig.~\ref{fig:motivation:initial:scaling}. The  level-based implementation saturates at a lower level, 
although both variants attain the same minimum code balance of $B_C=6\,\BF$ (Fig.~\ref{fig:motivation:initial:traffic}). 
The characteristic behavior is the same for the LB MPK with $p_m=4$ (squares): 
In line with the expectation, our method reduces the data traffic by a factor of approximately four ($B_{C,m} \approx 1.5\,\BF$) but it fails to improve performance at the full socket level. It even falls behind the baseline MPK for larger core counts. 
Further analysis reveals a 1.6$\times$ increase in retired instructions\footnote{using the event \texttt{INSTR\_RETIRED\_ANY} in \likwidperfctr} for LB MPK ($p_m=4$) compared to the baseline approach. These instructions are executed in the spin-waiting loop of OpenMP barriers~\cite{perf_pattern},
indicating that the synchronization between threads (performed after each computation of a level) is a potential bottleneck. 
An analysis of the level structure of the \texttt{pwtk} matrix confirms the relevance of synchronization cost as the average level size is approximately $850$ rows only. At an average of $53$ nonzeros per row, the workload of a level is just too low to ignore the synchronization cost, which increases with thread count and may reach a few thousand cycles at a full socket.\footnote{For the full CLX socket (ignoring hyper-threading) and the software environment used, a minimum barrier cost of 2,900 cycles was measured by direct barrier benchmarking.}  

The \texttt{Flan\_1565} matrix shows an opposite characteristic. 
The performance of LB MPK with $p_m=1$ is in line with the baseline approach, and the level blocking with $p_m=4$ achieves a performance improvement of 1.2$\times$ (see Fig.~\ref{fig:motivation:initial:scaling:Flan}). 
The moderate speedup of LB MPK is reflected in Fig.~\ref{fig:motivation:initial:traffic:Flan} by its rather high (measured) code balance of approximately $4\,\BF$, indicating that level-blocking is not very cache efficient in this case. 
The matrix level structure plays a decisive role here as there is a rather small number of levels, some of them being large. 
Already one of these large levels, which may contain up to $20,000$ rows (with about $75$ nonzeros per row) has a size of roughly $18\,\MB$, which is more than half of the L3 cache size of the CPU. 
Moreover, the small number of levels in combination with imbalanced level sizes may cause the irregular performance scaling of LB MPK ($p_m=4$) in Fig.~\ref{fig:motivation:initial:scaling:Flan}.
 

In the following three sections we describe three optimizations of the LB MPK, which are motivated by the performance shortcomings identified above. 
The first two are targeted at reducing the synchronization cost by forming larger levels (``level groups'') and substituting the expensive barrier by point-to-point synchronization.
The third optimization improves performance on matrices with dominant, bulky levels by recursively splitting these up (``recursion'') to improve cache efficiency.


\begin{figure*}[tbp]
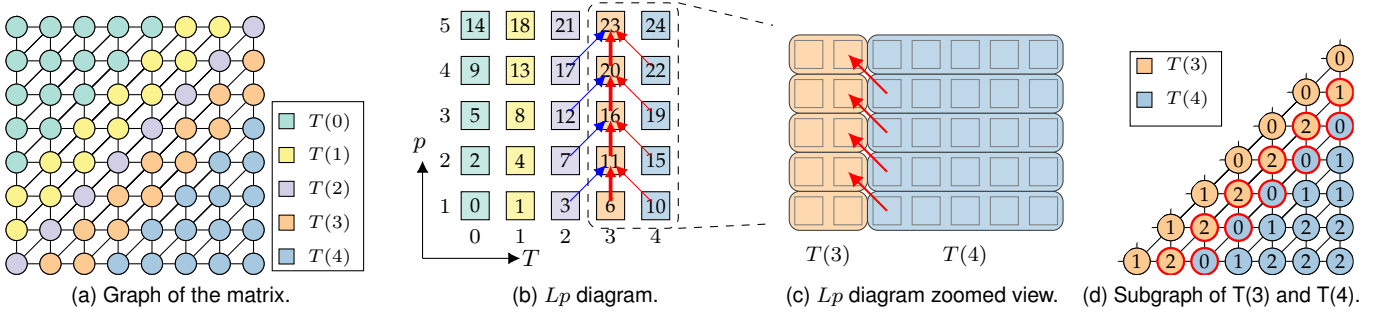

	\centering
	\subfloat[Graph of the matrix.]{%
	\input{plots/tikz/2d-7pt/graph_2d_7pt_permuted_w_levels_lg.tex}%
\label{fig:graph:lg}}
	\hfill
	\subfloat[$Lp$ diagram.]{%
	\input{plots/tikz/LP_graph/lp_group.tex}%
\label{fig:lp_graph:lg}}
	\hspace{-0.8em}
	\subfloat[$Lp$ diagram zoomed view.] {%
	\input{plots/tikz/LP_graph/lp_group_p2p_zoomed.tex}%
\label{fig:lg_graph:lg:zoomed}}
	\hfill
	\subfloat[Subgraph of T(3) and T(4).]{%
	\input{plots/tikz/2d-7pt/graph_2d_7pt_permuted_w_levels_lg_zoomed.tex}%
\label{fig:graph:lg:zoomed}}
	\caption{(a) Levels as in Fig.~\ref{fig:stencil_graph_permuted} being consolidated to five level groups ($T(0)$ -- $T(4)$). 
	(b) $Lp$ diagram corresponding to the level groups and
	the execution order of each level group at different power stages. 
	The bold red arrow (vertical) corresponds to the dependency with all the levels of the same level group $T(i)$ at the previous power stage $p-1$,  
	and the  slanted red arrow corresponds to the dependency  with the lowest-indexed level of next level group $T(i+1)$ at the previous power stage.
	The blue arrow corresponds to a dependency that is automatically fulfilled by the execution order.
	(c) Zoomed-in view of the $T(3)$ and $T(4)$ level groups in the $Lp$ diagram. The levels within the level group are seen as square nodes and the dependency between levels in $T(i)$ and $T(i+1)$ are clearly visible.
	The subgraph corresponding to the zoomed region is shown in (d). 
	The vertices drawn with red circles correspond to the two boundary levels between which synchronization in southeast direction has to be established.
	The numbers on the vertices represent the id of the thread (tid) working on that
	vertex.
\label{fig:level_group}}
\end{figure*}

\subsection{Level groups (LG)}
\label{sec:levelgoup}
The formation of larger levels follows the idea presented in~\cite{RACE}:
Successive levels are aggregated into so-called \emph{level groups}. 
This allows our LB MPK to operate on these level groups instead of the original levels.  
%
Figure~\ref{fig:graph:lg} shows the fifteen levels of Fig.~\ref{fig:stencil_graph_permuted} being clustered into five level groups $T(0)$--$T(4)$ ($T(i)$ denotes $i$-th level group). 
The $Lp$ diagram can easily be adapted by replacing the levels by the level groups on the $x$-axis (see Fig.~\ref{fig:lp_graph:lg}).\footnote{For the sake of uniformity we keep the name ``$Lp$'' for the diagram instead of ``$Tp$,'' although here we plot level groups ($T$) instead of levels ($L$) on the $x$-axis.}
Still, the same parallelogram-style blocking can be applied by traversing the level groups using the same rules as for the levels. 
Parallel execution is performed within a level group, and all threads synchronize after the computation of each group.
This strategy satisfies the neighborhood dependencies between levels as required by the LB MPK.


The cache reuse requirements of the LB MPK impose strict limits on the size of the level groups. 
As discussed in Sec.~\ref{sec:level_blocking_idea}, $p_m+1$ neighboring level groups have to be kept in cache. 
Therefore, if we assume neighboring level groups to be of similar
size, the following criterion has to be satisfied by the $i$-th level group $T(i)$:
\begin{equation}
(p_m+1)\times\NNZ(T(i))\times 12\,\bytes < f  C,
\label{eq:lc}
\end{equation}
where $\NNZ(T(i))$ is the number of nonzeros in $T(i)$, $C$ is a parameter representing the available cache size (in \bytes), and $f$ is a safety factor.
The cache size parameter is typically chosen to be less than or equal to the physical size of the cache(s) targeted for level blocking. The safety factor ($f=0.5$ in this work) accounts for extra traffic from other data structures and inefficiencies of the cache replacement policies.
The left part of inequality~\eqref{eq:lc} is the total memory traffic generated 
by accessing $p_m$+1 level groups (assuming 12\,\bytes\ per nonzero entry of the matrix, see Sec.~\ref{sec:baseline_MPK}), and the right part is the effective cache size.
If \eqref{eq:lc} is satisfied then level group $T(i)$
can be reused from cache for $p_m>1$; otherwise, at least parts of it must be loaded from main memory. 

Inequality~\eqref{eq:lc} is crucial to the construction process of the level groups.
We form the first level group $T(0)$ by accumulating
levels $L(0)$ \ldots\ $L(j)$  up to the largest $j$ for which  $\NNZ(L(0)) +  \dots + \NNZ(L(j)) = \NNZ(T(0))$ satisfies~\eqref{eq:lc}. 
The same procedure is repeated starting from 
level $L(j+1)$ to find $T(1)$, and successively forming the other level groups.
It can be seen from Fig.~\ref{fig:graph:lg} that this procedure
creates level groups with almost equal numbers of nonzero elements. 
In regions where levels contain fewer nonzeros per level,
more levels are aggregated (see $T(0)$ in Fig.~\ref{fig:graph:lg}) while in regions with bulkier levels, 
even a single level can form a level group (see $T(2)$ in Fig.~\ref{fig:graph:lg}).
As a result, the number of level groups is typically much smaller than the number of levels.

\begin{figure}[tbp]
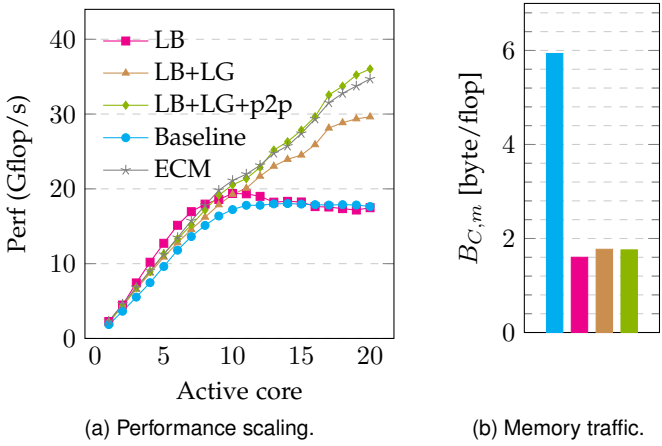

	\centering
	\subfloat[Performance scaling.]{%
	\input{plots/tikz/mtxPower_results/motivation/plot_opt_combined_pwtk.tex}%
\label{fig:motivation:opt:scaling}}
	\hfill
	\subfloat[Memory traffic.]{%
	\input{plots/tikz/mtxPower_results/motivation/plot_opt_combined_pwtk_traffic.tex}%
\label{fig:motivation:opt:traffic}}
	\caption{(a) Performance improvement of LB MPK using level group (LG) optimizations and point-to-point synchronization (p2p) for the \texttt{pwtk} matrix with $p_m=4$ on CLX.
	(b) Memory traffic of the four variants shown in (a). }\label{fig:motivation:opt}
\end{figure}

In the level-group-based scheme, synchronization only happens after the computation of each level group, which greatly diminishes the impact of barriers in case of LB MPK for the  \texttt{pwtk} matrix:
The performance of the LB MPK for $p_m=4$ (LB+LG; triangles in Fig.~\ref{fig:motivation:opt:scaling}) improves on a full socket by $1.6\times$ compared to the baseline MPK approach. 
At the same time, we only encounter a minor increase in the measured code balance (see Fig.~\ref{fig:motivation:opt:traffic}) since the condition~\eqref{eq:lc} limits the size of the level groups.  
Also the overhead from extra instructions reduces from 60\% for the naive LB MPK version to only 7\%.
The cache size parameter $C=35\,\MB$ has been set to the aggregate physical size of L3 and L2 caches of CLX. 

\subsection{Point-to-point (p2p) synchronization}
\label{sec:p2p}

The concept of level groups allows us to relax the lockstep-like synchronization by eliminating the OpenMP barrier after computation of a level group. 
The parallel LB MPK must ensure that the computations on the following levels and level groups are completed before the computation of power $p$ for a given level group $T(i)$: (A) the same level group $T(i)$ with previous power $p-1$ (bottom neighbor in $Lp$ diagram), (B) the highest-indexed (rightmost) level of $T(i-1)$ with power $p-1$ (southwest neighbor in $Lp$ diagram), and (C) the lowest-indexed (leftmost) level of $T(i+1)$ with power $p-1$ (southeast neighbor in $Lp$ diagram).

Note that the most stringent condition (A) can be enforced without a global barrier synchronization since it is only relevant when a level group $T(i)$ is visited again for computing the next power on it, which happens after a full diagonal traversal.
We thus implemented a customized locking mechanism (see below for details), which allows threads to spread out over a full diagonal of the $Lp$ diagram. They only need to check if all threads have finished computing the previous power of the current level group.
Due to the diagonal traversal scheme of the $Lp$ diagram, condition (A) implies condition (B) as the southwest neighbor of a level group  is always visited before its bottom neighbor (see numbering of execution order in ~\ref{fig:lp_graph:lg}). 
Finally, a similar mechanism is required to ensure condition (C). 
Here, only the completion of the relevant boundary level of the southeast neighbor has to be ensured (see Fig.~\ref{fig:lg_graph:lg:zoomed}).
As vertices are statically assigned to threads, this boundary level is typically calculated by a single or only a few threads, further relaxing the synchronization demands between all threads. This can be observed in Fig.~\ref{fig:graph:lg:zoomed}, where only the first thread (tid=0) is involved in the computation of the relevant boundary level of $T(4)$.

 The locking mechanism is implemented as follows: An array of locks is defined in the initialization phase, which allows us to control the above dependencies.
The threads use \texttt{omp atomic} and spin-waiting loops to set and test the locks. 
Figure~\ref{fig:motivation:opt:scaling} shows the performance scaling of this implementation (LB+LG+p2p; diamonds) in comparison to the other variants; it yields a performance boost of $1.2\times$ over the version with level groups and barrier synchronization (LB+LG).
A part of this speedup comes from the reduced synchronization cost. The 
rest is due to the relaxation of lock step synchronization that
allows for overlap between memory and cache transfers, i.e.,
some threads can work on the memory-bound phase ($p=1$) while the rest work on a cache-bound phase ($p>1$).
The optimization thus brings us close to our phenomenological ECM model
(stars in  Fig.~\ref{fig:motivation:opt:scaling}) and results in a $2\times$ speedup over the baseline approach.
	Note that as the sizes of level groups change, traffic within inner cache levels will also change. Since the ECM model uses this data traffic as input, 
	it results in slightly different models when sizes of level groups change.
	This can be observed for example by comparing Fig.~\ref{fig:motivation:initial:scaling} and Fig.~\ref{fig:motivation:opt:scaling}.

\begin{figure}[tbp]
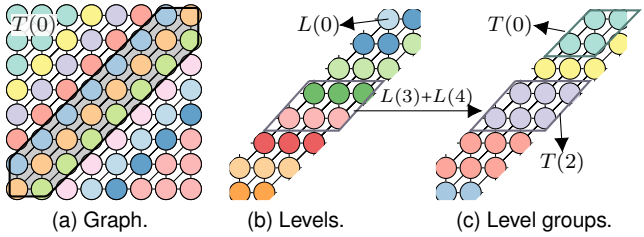

	\centering
	\subfloat[Graph.]{%
	\input{plots/tikz/2d-7pt/graph_2d_7pt_permuted_w_levels_lg_rec_stage0.tex}%
\label{fig:graph:rec:stage0}}
	\hfil
	\subfloat[Levels.\textcolor{white}{please align}]{%
	\input{plots/tikz/2d-7pt/graph_2d_7pt_permuted_w_levels_rec_stage1.tex}%
\label{fig:graph:rec:stage1}}
	\hspace{-3em}
	\subfloat[Level groups.]{%
	\input{plots/tikz/2d-7pt/graph_2d_7pt_permuted_w_levels_rec_stage1_lg.tex}%
\label{fig:graph:rec:stage1:lg}}
	\caption{(a) Level groups in the graph. The shaded subgraph shows 
		the level groups with more than six rows, where recursive treatment is applied.
		(b) BFS levels within the subgraph. (c) Level groups formed from the levels within the subgraph. 
 		\label{fig:rec}}
\end{figure}

\subsection{Recursion}
\label{sec:rec}

The negative impact of bulky levels (which do not satisfy~\eqref{eq:lc}) on main memory traffic for the  LB MPK approach (see Fig.~\ref{fig:motivation:initial:traffic:Flan}) has been identified and discussed for the \texttt{Flan\_1565} matrix in Sec.~\ref{sec:level_blocking_idea}.
In the RACE coloring scheme~\cite{RACE}, a recursive approach has been presented to generate higher levels of parallelism within bulky levels. 
The same method can be used in our context to successively generate new levels or level groups of reduced size until they fit into cache. 
The idea is to apply the LB MPK presented so far to the subgraph defined by a single level or a set of consecutive levels. 
As a result, a new set of smaller levels is generated for this subgraph. 
If some of the new levels still violate~\eqref{eq:lc}, the procedure is applied again to the new subgraph defined by these levels.  
This procedure can be continued until all levels fit into a cache.



\begin{figure}[tbp]
	\subfloat[Without recursion ($s_m=0$).]{\includegraphics{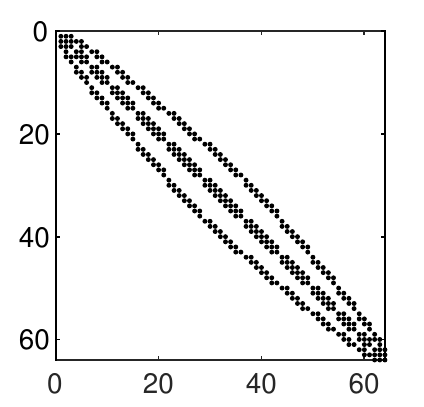}
		\label{fig:matrix_sparsity:wo_recursion}}
	\hfill
	\subfloat[With recursion($s_m=1$).]{\includegraphics{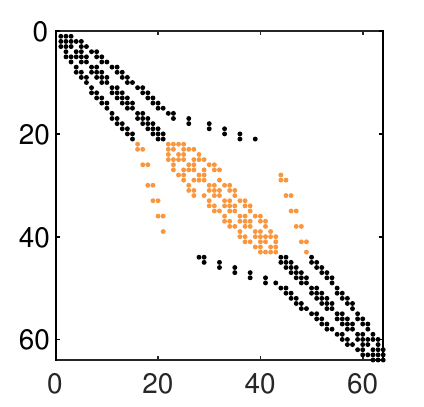}\label{fig:matrix_sparsity:w_recursion}}
	\caption{
		Sparsity pattern of the stencil example matrix without (a) and with (b) recursion. The entries of submatrix where recursion is applied is shown with orange color in (b). \label{fig:matrix_sparsity}}
\end{figure}

We start by locating (consecutive) levels that do not fit in a cache and isolate the subgraph formed by these levels.  
BFS is applied first to this subgraph, and then a set of level groups is formed from these BFS levels.
The resulting level groups are typically smaller than the previous ones as neighboring vertices outside the subgraph do not need to be considered.
Figure~\ref{fig:rec} illustrates this procedure for our stencil example and a hypothetical cache which
satisfies~\eqref{eq:lc} for level groups $T(i)$ containing no more than six vertices.
We find that the three bulkier level groups (containing one level each) $T(4)$ -- $T(6)$ do not satisfy the condition. 
The subgraph  induced by these three levels is formed (shaded with gray background in Fig.~\ref{fig:graph:rec:stage0}),  and we identify the eight BFS levels of this subgraph (Fig.~\ref{fig:graph:rec:stage1}).
Following the discussion in Sec.~\ref{sec:levelgoup}, the level groups of the subgraph are constructed (Fig.~\ref{fig:graph:rec:stage1:lg}). 
They are now small enough to satisfy~\eqref{eq:lc} and the process stops.


In general, the procedure can be applied recursively until the level groups satisfy~\eqref{eq:lc} or
a user-specified maximum recursion stage $s_m$ is reached, where $s_m=0$ is the case without any recursion. 
In the following, $s$ ($\le s_m$) denotes the current recursion stage.
The maximum recursion stage should, however, be limited as applying the recursion step leads to loss of data locality at the boundaries of the subgraph. 
This happens because the subgraphs are permuted (BFS) without taking into account the neighbors outside the subgraph.
Figure~\ref{fig:matrix_sparsity} demonstrates this effect by comparing the matrix structure of our stencil example without recursion ($s_m=0$) and with one recursion step ($s_m=1$) applied to the inner levels. 
The matrix bandwidth increases for the boundary elements of the subgraph because of the mismatch of the vertex numberings outside and inside the subgraph. 
While access to the matrix elements remains linear, the more irregular accesses to the right-hand side vector may impact the overall MPK performance.
Note that the graphical representation in Fig.~\ref{fig:matrix_sparsity:w_recursion} exaggerates this effect, since in our toy problem the subgraph represents a substantial fraction of the full problem.
%
The performance influence of the maximum recursion stage $s_m$ is discussed later in Sec.~\ref{sec:s_m_influence}.


\begin{figure}[tbp]
	\centering
	%
	\input{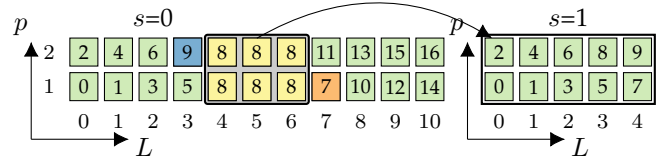}%

	\caption{The $Lp$ diagram for $p_m=2$. \emph{Left}: $Lp$ diagram of the $s=0$ recursion stage ($Lp^0$), which contains level groups of the entire graph seen in Fig.~\ref{fig:graph:rec:stage0}.
	The level groups selected for recursion are highlighted. 
	\emph{Right}: $Lp$  diagram at $s=1$ ($Lp^1$),  which consists of the level groups shown in Fig.~\ref{fig:graph:rec:stage1:lg}.
	The execution order of the $Lp$ graph is shown with numbers.
\label{fig:lp:rec:p2}
	}	
\end{figure}

\begin{figure}
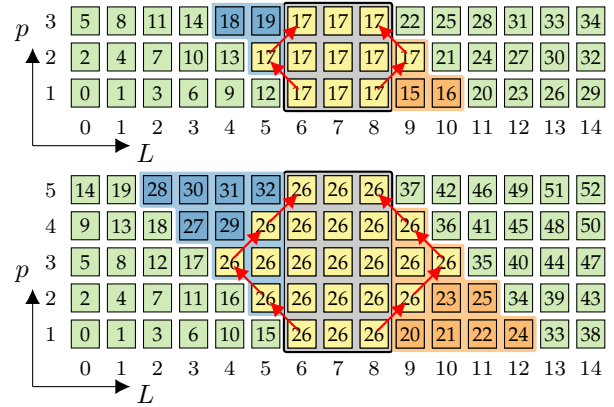

 	\centering
	%
	\input{plots/tikz/LP_graph/1_subgraph_in_one_rec_stage/lp_recursion_p3_diamond.tex}%

	%
	\input{plots/tikz/LP_graph/1_subgraph_in_one_rec_stage/lp_recursion_p5_diamond.tex}%

	\caption{$Lp^0$ diagrams with $p_m=3$ (above) and $p_m=5$ (below) corresponding to an arbitrary graph where recursion has to be applied to level groups $T(6)$--$T(8)$ forming $Lp^1$ (not shown).
	The red arrows show the longest input (output) dependency from (to) the boundary points of recursive region.
	\label{fig:lp:rec:diamond}}	
\end{figure}

As each subgraph (formed from consecutive levels) of a recursion stage creates its own level groups, we construct $Lp$ diagrams for each subgraph, i.e., $Lp^s$ represents the $Lp$ diagrams of recursion stage $s$.
Figure~\ref{fig:lp:rec:p2} shows the two $Lp$ diagrams of the stencil example for $p_m=2$: $Lp^0$ representing $s=0$ on the full graph (Fig.~\ref{fig:graph:rec:stage0}), and $Lp^1$ after the first recursion stage of the subgraph corresponding to level groups in Fig.~\ref{fig:graph:rec:stage1:lg}.
Note that the numbering of the execution order is local to each $Lp^s$ diagram. 
All level groups of a subgraph of $Lp^s$ to which recursion is applied have the same execution order in $Lp^s$ (e.g., the subgraph related to $T(4)$--$T(6)$ in $Lp^0$ is executed in step $8$ of $Lp^0$ in Fig.~\ref{fig:lp:rec:p2}). 
The actual execution order of the vertices in this subgraph is determined by $Lp^{s+1}$ (see $Lp^1$ in Fig.~\ref{fig:lp:rec:p2}). 
In general, the actual execution of a given vertex is determined by the $Lp$
diagram associated with the highest recursion stage of the vertex.
Of course the actual execution order in the $Lp^s$ diagrams still needs to maintain the data dependencies of the LB MPK. 
With $p_m=2$ as used in Fig.~\ref{fig:lp:rec:p2} we can still maintain our diagonal-type execution order within the graphs:  $T(7)$ of $Lp^0$ is updated to $p=1$ at step $7$. $Lp^1$ is calculated as step $8$ of $Lp^0$. In step $9$ of $Lp^0$,  $T(3)$  is updated to $p=2$.

For $p_m>2$, the dependency relations between execution order of $Lp^s$ and $Lp^{s+1}$ are more complicated. This is depicted in Figure~\ref{fig:lp:rec:diamond}, where $Lp^0$ with  $p_m=3,5$ is shown for 15 level groups and $T(6)$--$T(8)$ form the subgraph on which $Lp^1$ is built. 
Actually, all nodes in the parallelogram formed by the diagonals in $Lp^s$ ($Lp^0$ in our example) that cross the subgraph to be refined have dependency relations to the vertices in this subgraph. 
Within the parallelogram, there are three different types of dependencies for the nodes to be computed at $Lp^s$ ($Lp^0$ in Figure~\ref{fig:lp:rec:diamond}) and which are not in the subgraph to be refined:
(i) Nodes which provide input only to $Lp^{s+1}$ and which need to be calculated before $Lp^{s+1}$ (orange color in Figure~\ref{fig:lp:rec:diamond}),
(ii) nodes which have only an output dependency on $Lp^{s+1}$ and need to be calculated after $Lp^{s+1}$ (blue color in Figure~\ref{fig:lp:rec:diamond}), 
(iii) nodes within the ``diamond'' embedded in the parallelogram, which have input and output dependencies related to the computations in $Lp^{s+1}$ and need to be calculated in coordination with $Lp^{s+1}$.  
All nodes within the ``diamond'' thus have the same execution order in $Lp^s$, and the calculation of $Lp^{s+1}$ also involves computation of level groups outside the subgraph to be refined. 
This diamond-type execution structure is well known from diamond tiling~\cite{diamond_blocking} applied to stencils.

Note that this recursive refinement approach is not limited to a single subgraph of a given $Lp^s$. However, if multiple subgraphs need to be refined, the parallelograms formed by these subgraphs must not overlap.

\begin{figure}[tbp]
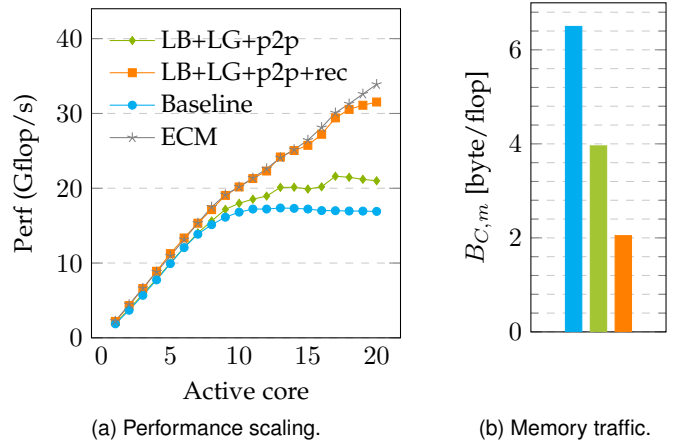

	\centering
	\subfloat[Performance scaling.]{%
	\input{plots/tikz/mtxPower_results/motivation/plot_opt_combined_Flan_1565.tex}%
\label{fig:motivation:rec:scaling}}
	\hfill
	\subfloat[Memory traffic.]{%
	\input{plots/tikz/mtxPower_results/motivation/plot_opt_combined_Flan_1565_traffic.tex}%
\label{fig:motivation:rec:traffic}}
	\caption{(a) Performance improvement of LB MPK using recursion (squares) compared to the one without recursion (diamonds) for the \texttt{Flan\_1565} matrix with $p_m=4$ on one socket of CLX.
		Both versions use level groups and p2p optimizations.
		The performance of the baseline approach as well as the ECM model is also shown for reference. (b) Measured memory traffic of the three variants on the left.
	}\label{fig:motivation:rec}
\end{figure}
The impact of the presented recursion scheme on the performance of the LB MPK method for  the \texttt{Flan\_1565} matrix with $p_m=4$ is shown in Fig.~\ref{fig:motivation:rec:scaling}.
We used a cache size parameter $C=45\,\MB$ for LB MPK methods and set $s_m=4$ for the case with recursion (squares). In this setting, the $Lp^0$ diagram has three subgraphs to which recursive treatment is applied.
Via improved cache reuse, the recursion improves the full-socket performance by a factor of almost $1.4\times$ compared to the version without recursion.
This comes with a corresponding reduction of almost $2\times$ in main memory data traffic (Fig.~\ref{fig:motivation:rec:traffic}). 
Compared to the baseline MPK approach, we achieve an overall reduction of main memory traffic by  $3.2\times$ and an increase in performance by $1.8\times$ on a full socket of CLX. 
These numbers and the (close to) linear scaling of our method indicate that main memory access is no longer the performance bottleneck.

\subsection{RACE}
The LB MPK algorithm including all optimizations discussed above has been implemented in the RACE library (code available at~\cite{RACE-git}).
In the following we therefore refer to our LB MPK implementation as ``RACE MPK.''
The library supports both preprocessing and execution phases of the LB MPK.
For preprocessing, RACE requires the matrix, highest power $p_m$, cache size $C$, and maximum recursion stage $s_m$ as input  and returns the permutation vector as output. 
The user then has to pass the permuted matrix and a call-back function 
to RACE for execution. 
RACE will execute the call-back function in parallel (using OpenMP threading) according to the 
internally created \texttt{level\_ptr} and $Lp$ diagrams.
\begin{figure*}[tbp]
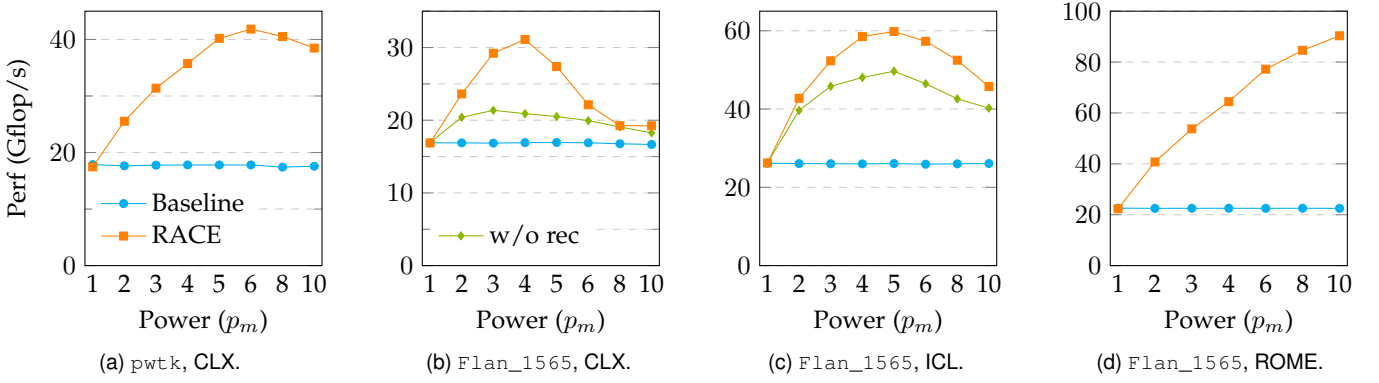

	\centering
	\subfloat[\texttt{pwtk}, CLX.]{%
	\input{plots/tikz/param_study/plot_p_pwtk.tex}%
\label{fig:param_study:p:pwtk}}
	\hfill
	\subfloat[\texttt{Flan\_1565}, CLX.]{%
	\input{plots/tikz/param_study/plot_p_Flan_1565.tex}%
\label{fig:param_study:p:Flan}}
	\hfill
	\subfloat[\texttt{Flan\_1565}, ICL.]{%
	\input{plots/tikz/param_study/plot_p_Flan_1565_ICL.tex}%
\label{fig:param_study:p:Flan:ICL}}
	\hfill
	\subfloat[\texttt{Flan\_1565}, ROME.]{%
	\input{plots/tikz/param_study/plot_p_Flan_1565_ROME.tex}%
\label{fig:param_study:p:Flan:ROME}}
	\caption{Performance as a function of maximum power $p_m$
	  for RACE and the baseline implementation of MPK.
	  For cases where recursion yields a speedup,
	  we also plot the performance of RACE without recursion (in green) for comparison.
	}\label{fig:param_study:p}
\end{figure*}


\section{Parameter study}
\label{sec:param_study}
Our RACE MPK as introduced in the previous section has three input parameters: the 
maximum power $p_m$, the cache size $C$, and the maximum recursion stage $s_m$.
In this section we discuss the qualitative impact of these parameters on the performance of RACE MPK.

\subsection{Influence of \boldmath\textbf{$p_m$}}
\label{sec:power_influence}
Ideally, RACE MPK requires to access main memory for each level group exactly once at $p=1$.  
The remaining $p_m-1$ accesses can potentially be served from the cache(s) (see Figs.~\ref{fig:motivation:opt:traffic} and~\ref{fig:motivation:rec:traffic}).
As a consequence, cache utilization and performance should increase with $p_m$. 
However, as $p_m$ gets larger, the number of level groups grows and their size must reduce as condition~\eqref{eq:lc} has to be fulfilled, which results in higher synchronization cost.
These opposing effects result in a typical performance pattern as shown in 
Fig.~\ref{fig:param_study:p:pwtk} for the \texttt{pwtk} matrix on CLX.
Initially the performance increases almost linearly with $p_m$ but starts to drop gradually at larger $p_m$ ($\approx 6$--$8$ in our example).
For matrices that require recursion, the performance drop is more prominent and occurs at a lower $p_m$ as shown in Fig.~\ref{fig:param_study:p:Flan}  for the \texttt{Flan\_1565} matrix on CLX.
The additional overhead at the boundaries of the recursively refined level groups (see discussion in Sec.~\ref{sec:rec}) add another performance penalty.
Of course, the $p_m$ value at which performance starts to decrease 
depends on the matrix and the cache size.
This can be observed by comparing the performance of \texttt{Flan\_1565} 
on the three architectures (Figs.~\ref{fig:param_study:p:Flan}--\ref{fig:param_study:p:Flan:ROME}). 
On ROME (Fig.~\ref{fig:param_study:p:Flan:ROME}) with its large last-level cache, the matrix does not require
recursion at all and the performance increases up to $p_m=10$,
where the RACE MPK achieves a speedup of 4$\times$ compared to the baseline MPK.
The ICL (Fig.~\ref{fig:param_study:p:Flan:ICL}) and CLX (Fig.~\ref{fig:param_study:p:Flan}) CPUs need recursion to achieve best performance. 
The maximum performance is attained at $p_m$ values of 5 and 4, resulting in speedups of 2.3$\times$ and 1.8$\times$ with respect to the MPK baseline on these two architectures. Note that performance improvements decrease with decreasing cache sizes.

For applications computing $A^{k}x$ using RACE MPK, the best strategy is to identify the optimal $p_m$ value $p_m^{\mathrm{opt}}$ and perform the $A^{p_m^{\mathrm{opt}}}x$ computations multiple times (if $k$$>$$p_m^{\mathrm{opt}}$) until the power $k$ is reached.
If $k$ is not a multiple of $p_m^{\mathrm{opt}}$, the remainder computations can be done using MPK kernels with $p_m < p_m^{\mathrm{opt}}$.
\begin{figure}
	\centering
	\subfloat[\texttt{pwtk}.]{\includegraphics[width=0.5\linewidth]{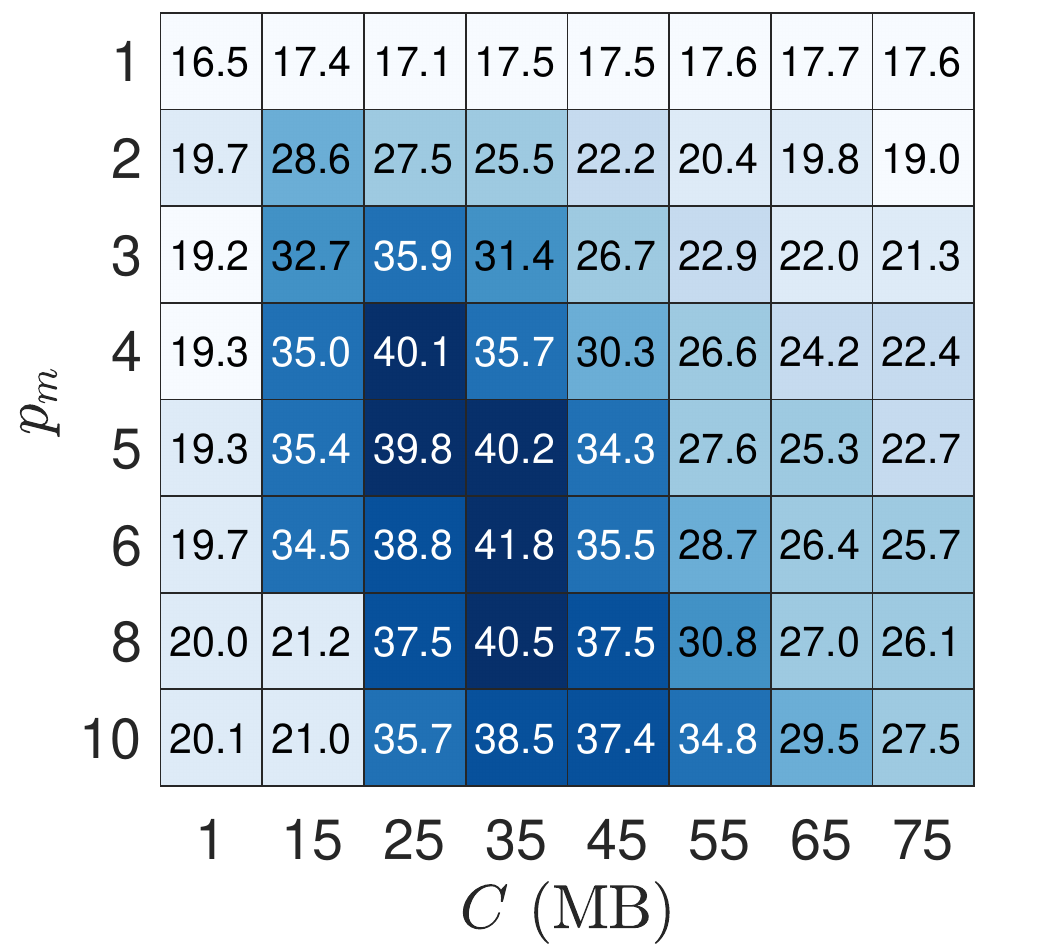}\label{fig:param_study:cs:pwtk}}
	\hfill
	\subfloat[\texttt{Flan\_1565}.]{\includegraphics[width=0.5\linewidth]{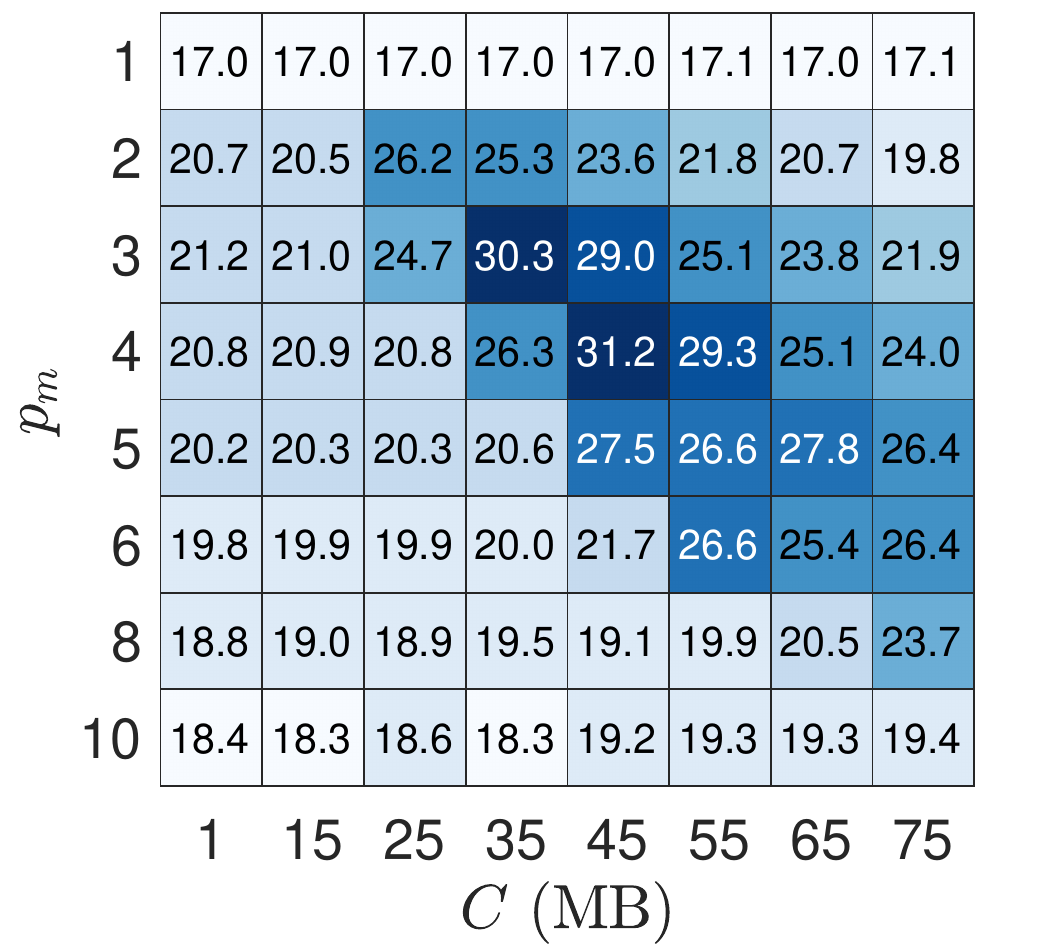}\label{fig:param_study:cs:Flan}}
	\caption{Influence of cache size $C$ and power $p_m$ on performance (in \GFS) of the RACE MPK using all cores of CLX.}\label{fig:param_study:cs}
\end{figure}

\subsection{Influence of \boldmath\textbf{$C$}}
The interaction of cache size $C$ and highest power $p_m$ is shown as a heatmap in Fig.~\ref{fig:param_study:cs} for the \texttt{pwtk} and \texttt{Flan\_1565} matrices on CLX.
The optimal $C$ value is between 25 and 45\,\MB\ irrespective of $p_m$ and the matrix. 
This is in good qualitative agreement with the aggregate size of the L3 (27.5\,\MiB, victim) and L2 cache (20\,\MiB) of CLX.
Of course, the RACE MPK method works best when blocking for the biggest available cache.
Smaller $C$ values lead to 
smaller level groups (see~\eqref{eq:lc}) and therefore higher synchronization and recursion overheads.
On the other hand, $C$ values bigger than the total cache size will obviously provoke cache misses.
\begin{figure}
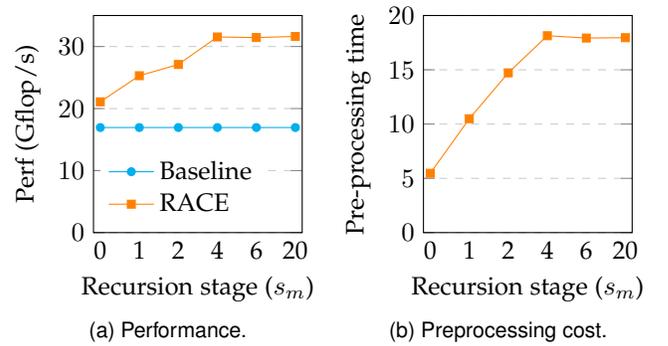

	\centering
	\subfloat[Performance.]{%
	\input{plots/tikz/param_study/plot_recStage_Flan_1565.tex}%
\label{fig:param_study:rec_stage:perf}}
	\hfil
	\subfloat[Preprocessing cost.]{%
	\input{plots/tikz/param_study/plot_recStage_Flan_1565_pre_time.tex}%
\label{fig:param_study:rec_stage:pre_time}}
	\caption{(a) Performance influence of maximum recursive stage $s_m$ on the performance of the \texttt{Flan\_1565} matrix with $p_m=4$ and $C=35$\,\MB\ on one socket of CLX. (b) Corresponding preprocessing cost of RACE in equivalent number of SpMVs.\label{fig:param_study:rec_stage}}
\end{figure}

\subsection{Influence of \boldmath\textbf{$s_m$}}
\label{sec:s_m_influence}
For matrices that require recursion to fulfill~\eqref{eq:lc},
the maximum recursion depth $s_m$ may stop the recursion procedure even if the condition is still violated for some level groups. 
Figure~\ref{fig:param_study:rec_stage:perf} depicts the performance behavior of the
\texttt{Flan\_1565} matrix with $p_m=4$ on CLX as a function of $s_m$.
Initially, the performance increases with $s_m$ as the level groups become smaller. When~\eqref{eq:lc} is fulfilled at $s_m=4$ for all level groups, performance saturates. 
Note that increasing $s_m$ does not always have the positive performance effect as observed for \texttt{Flan\_1565}. 
The overhead at the boundaries of the refined subgraphs may overcompensate the gains of increased cache efficiency.  
For example, in case of the \texttt{RM07R} matrix on ICL (not shown in plots) with $p_m=3$ ($=p_m^{\mathrm{opt}}$) it was found that $s_m=0$ (no recursion) achieves  1.2$\times$ better performance than $s_m=13$, where all the level groups fit in cache.
Of course, the optimal value of $s_m$ is determined by an intricate interplay of cache properties and matrix properties and thus cannot be found analytically. Typically, recursion should only be applied if condition~\eqref{eq:lc} cannot be fulfilled with $s_m=0$. In this scenario, recursion depths up to $s_m=15,\ldots,20$ should be scanned for best performance.  
	

The preprocessing cost increases with $s_m$ as levels have to be found  for recursive subgraphs.
This can be seen in Fig.~\ref{fig:param_study:rec_stage:pre_time} for the \texttt{Flan\_1565} matrix, where the
preprocessing cost (shown in equivalent SpMVs) increases with $s_m$ up to $s_m=4$. 
The construction of levels (BFS) dominates the preprocessing time.
The other parameters $p_m$ and $C$ do not have a considerable impact on preprocessing time as changing them does not require to generate new levels.

\section{Performance evaluation}
\label{sec:perf_eval}
In this section we investigate the performance of RACE MPK 
and compare it against the baseline MPK for 36 different sparse matrices commonly seen in literature.
The details of these matrices can be found in Table~\ref{tab:matrices}.

\subsection{Experimental setup}
\label{sec:mpk:exp_setup}
All matrices were stored in the CRS data storage format (see Sec.~\ref{sec:baseline_MPK}). 
Unless  specified otherwise we used all the cores on one CPU socket and one thread per core. 
To ensure vectorization of the kernels we used
\texttt{\#pragma simd vectorlength(VECLEN) reduction(+:tmp)} on the innermost loop of the SpMV (see Fig.~\ref{fig:SpMV_alg}).
The vector length (\texttt{VECLEN}) was specified explicitly and was chosen to be the maximum SIMD width of the hardware,
i.e., four on ROME and eight on ICL and CLX.

For both baseline and RACE, the matrices were preprocessed with RCM reordering using the Intel SpMP~\cite{SpMP} library if it improved the performance. 
The baseline method was parallelized using the \texttt{\#pragma omp parallel for schedule(static)} workshare construct along the outermost loop (over matrix rows).\footnote{Note that static scheduling was chosen as the benchmark matrices (see Table~\ref{tab:matrices}) did not have highly imbalanced row lengths.}
RACE is parallelized using OpenMP pragmas by manually assigning the vertices in each level group to the threads and implementing the point-to-point synchronization mechanisms discussed in Sec.~\ref{sec:p2p}.
The parameter space of RACE (see Sec.~\ref{sec:param_study}) was tuned in 
the following range: $p_m\in\{[\text{$1$:$1$:$3$}] \cup [\text{$4$:$2$:$16$}]\}$,\footnote{in the format [start value : increment : end value]} $C$ in the range of total cache (L3+L2) size of the hardware, and $s_m\in\{0,1,2,4,6,20\}$.
More specifically, the parameter space of $C$ (in \MB) is $[\text{$25$:$10$:$45$}]$ for CLX,
$[\text{$65$:$10$:$105$}]$ for ICL, and $[\text{$100$:$50$:$250$}]$ for ROME.
\begin{figure*}[tbp]
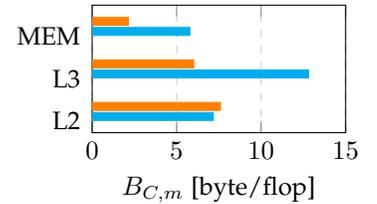
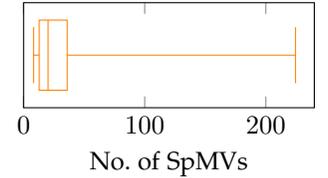
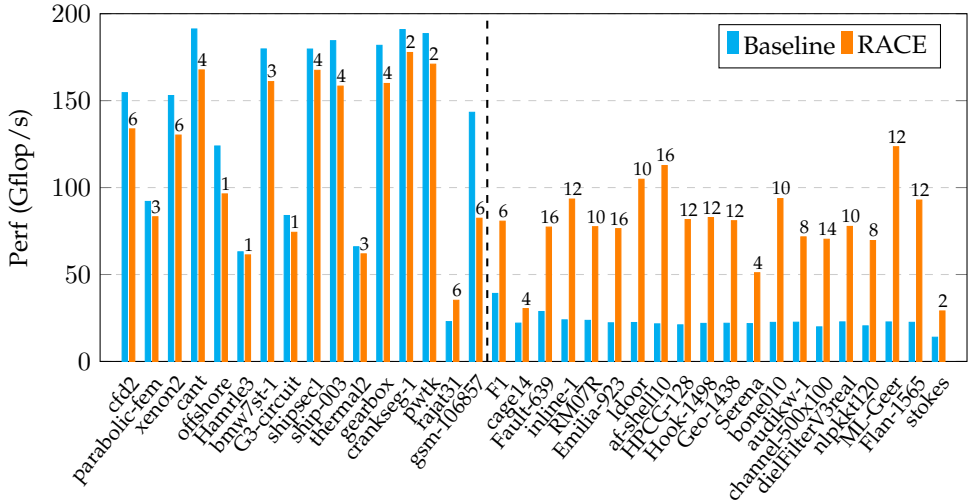
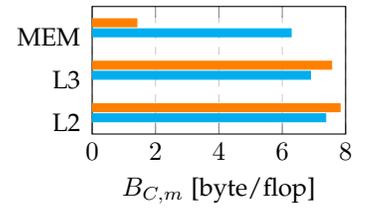
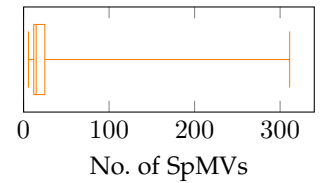

	\centering
	\begin{minipage}{0.72\textwidth}
	\subfloat[Performance, CLX.]{%
	\input{plots/tikz/mtxPower_results/plot_casclakesp2.tex}%

		\label{fig:mtxPower_results:clx:perf}}\quad
	\end{minipage}
	\begin{minipage}{0.272\textwidth}
	\subfloat[Data traffic, CLX.]{%
	\input{plots/tikz/mtxPower_results/plot_casclakesp2_traffic_rotated.tex}%

		\label{fig:mtxPower_results:clx:mem}}\\
	\subfloat[Preprocessing cost, CLX.]{%
	\input{plots/tikz/mtxPower_results/plot_casclakesp2_preTime_rotated.tex}%

	\label{fig:mtxPower_results:clx:pre_time}}
	\end{minipage}
	\hfil
	\begin{minipage}{0.72\textwidth}
	\subfloat[Performance, ICL.]{%
	\input{plots/tikz/mtxPower_results/plot_horeka.tex}%

		\label{fig:mtxPower_results:icl:perf}}
	\end{minipage}
	\begin{minipage}{0.272\textwidth}
	\subfloat[Data traffic, ICL.]{%
	\input{plots/tikz/mtxPower_results/plot_horeka_traffic_rotated.tex}%

	\label{fig:mtxPower_results:icl:mem}}\\
	\subfloat[Preprocessing cost, ICL.]{%
	\input{plots/tikz/mtxPower_results/plot_horeka_preTime_rotated.tex}%

	\label{fig:mtxPower_results:icl:pre_time}}
	\end{minipage}
	\hfil
	\begin{minipage}{0.72\textwidth}
	\subfloat[Performance, ROME.]{%
	\input{plots/tikz/mtxPower_results/plot_tg094.tex}%

		\label{fig:mtxPower_results:rome:perf}}
	\end{minipage}
	\begin{minipage}{0.272\textwidth}
	\subfloat[Data traffic, ROME.]{%
	\input{plots/tikz/mtxPower_results/plot_tg094_traffic_rotated.tex}%

		\label{fig:mtxPower_results:rome:mem}}\\
	\subfloat[Preprocessing cost, ROME.]{%
	\input{plots/tikz/mtxPower_results/plot_tg094_preTime_rotated.tex}%

		\label{fig:mtxPower_results:rome:pre_time}}
	\end{minipage}
	\caption{(a), (d), (g): Performance comparison between baseline and RACE MPK on CLX, ICL, and ROME, respectively. 
	The dashed line represents the total available cache size and the 
	numbers show the tuned $p_m$ values corresponding to the RACE performance. 
	(b),(e), (h): L2, L3, and memory code balance of 
	RACE MPK and baseline approach on the three architectures.
	The memory and cache data traffic shown is the average across all the in-memory matrices (i.e., to right of dashed line in the respective performance plot).
	(c), (f), (i): Statistics of the preprocessing cost of RACE MPK
	for all in-memory matrices. 
	The cost is shown as the number of SpMVs that can be executed in the given time.}
	\label{fig:mtxPower_results}
\end{figure*}

\subsection{Results}
\label{sec:mpk_results}
Figures~\ref{fig:mtxPower_results:clx:perf}, \ref{fig:mtxPower_results:icl:perf}, and \ref{fig:mtxPower_results:rome:perf} show the performance of baseline and RACE MPK on CLX, ICL, and ROME, respectively.
The matrices are ordered (left to right) according to increasing data-set size (number of nonzeros).
The vertical lines represent the total cache size of the respective hardware and thus categorize matrices into memory-resident (right of line) and  cache-resident (left of line) scenarios. 

For the smallest matrices, RACE does not usually show significant speedup over the baseline method
as these matrices comfortably fit in cache.  
However, as the working set approaches the cache size, RACE starts to develop clear performance advantages. 
On CLX and ICL, this effect is pronounced already for larger ``in-cache'' matrices, while for ROME the benefit of RACE MPK starts exactly at the boundary between cache- and memory-resident matrices.  
There are two main reasons for this: (i) The transition between L3 and main memory bandwidth on Intel architectures is gradual compared to AMD ROME (see Fig.~\ref{fig:bandwidth}), and
(ii) the L3 and L2 caches have almost similar sizes on both Intel architectures, and the blocking in RACE targets the combined L3 and L2 caches. 
Therefore,  for smaller matrices that fit into the L3 cache, RACE can reduce the L2 traffic compared to the baseline method. 
On the other hand, for ROME the L3 cache is considerably bigger than the L2 and hence the blocking is performed only in the L3 cache, thereby bearing no benefit for matrices fitting in the L3 cache.

For all memory-resident matrices RACE has a clear performance advantage on all architectures, achieving typical speedups of 2$\times$ to 5$\times$ compared to the baseline MPK.
This is correlated with the measurements of the average L2, L3, and main memory traffic shown in Figs.~\ref{fig:mtxPower_results:clx:mem}, \ref{fig:mtxPower_results:icl:mem}, and \ref{fig:mtxPower_results:rome:mem}.
Here the baseline MPK approach is close to the SpMV's minimum traffic limit of $6\,\BF$\footnote{The L3 traffic measurements using \likwidperfctr\ is double on CLX and ICL as the current version of \likwidperfctr\ cannot distinguish traffic between main memory and L2 cache with  L3 and L2 caches; see~\cite{likwid_L3_problem} for details.}, indicating
 the absence of caching of matrix elements. 
In most cases the baseline approach is also strongly memory bandwidth bound and thus performs close to the optimistic (memory-bound) \rl\ limit (i.e., ${b_\mathrm{Mem}}/{B_C}$)
of 19,  28, and 24\,\GFS\ on CLX, ICL, and ROME, respectively.
For RACE we find a memory traffic less than the minimum SpMV limit on all the three architectures due to caching of the matrix elements.
On CLX and ICL, even the L3 traffic reduces substantially as the large (aggregate) L2 cache contributes substantially to the blocking.
The reduced data traffic of RACE results in a  performance higher than the SpMV in-memory \rl\ limit and  the baseline approach.
Correlated with the reduction of main memory traffic, RACE achieves the highest speedups on ROME where we observe an average (maximum) performance gain of 3.5$\times$ (5.4$\times$).
On ICL and CLX, we observe an average speedup of almost 2$\times$ and 1.6$\times$, respectively, and a maximum speedup of 3$\times$ and 2.3$\times$.

The significantly higher performance (as well as speedup) of RACE on ROME compared to the Intel systems can be attributed to
its larger L3 cache and higher L3 bandwidth (see Fig.~\ref{fig:bandwidth}).
A larger L3 allows to  cache level groups for higher $p_m$ values (see~\eqref{eq:lc}).
This can be observed in the tuned $p_m$ values
annotated with numbers on top of the RACE performance bars.
We see that for the same matrices the $p_m$ values on ROME are higher than that of ICL and CLX. 
This allows for matrix elements to be cached longer on ROME and results in an average memory traffic reduction of 4.5$\times$ (see Fig.~\ref{fig:mtxPower_results:rome:mem}) compared to the baseline, while on ICL and CLX the reduction is 2.7$\times$ and 2.2$\times$, respectively.

\subsection{Preprocessing cost}
Now that the performance behavior of RACE is understood, we need to investigate its preprocessing overhead. 
The box plots in Figs.~\ref{fig:mtxPower_results:clx:pre_time}, \ref{fig:mtxPower_results:icl:pre_time}, and \ref{fig:mtxPower_results:rome:pre_time} 
show statistics of RACE's preprocessing cost for memory-resident matrices.
These cost is shown in equivalent number of SpMVs that 
can be executed during the time required for preprocessing.
In general, the cost reduces as the cache size of the architecture increases, i.e., on ROME the preprocessing time is well under the time of 30 SpMVs for most matrices while on ICL and CLX the equivalent SpMV invocations are 40 and 60, respectively.
This is due to larger cache sizes 
requiring fewer recursion stages ($s_m$), since the preprocessing cost increases with $s_m$ (see Fig.~\ref{fig:param_study:rec_stage:pre_time}).

Most of the preprocessing time ($>$95\%) is spent on determining
the levels using BFS. 
In RACE we use a parallel BFS implementation similar to
the top-down approach from~\cite{BFS_Scott},
where the parallelization is accomplished by distributing the 
vertices in a level (frontier) to different threads.
However, this method lacks sufficient parallelism if the number of vertices in a level is too small. 
This is the case with the \texttt{RM07R} matrix, which is the outlier in the preprocessing cost on all three architectures. 
Here, a lot of levels contain only one vertex and preprocessing is largely sequential. 

\section{Application: Chebyshev time propagation}
\label{sec:app_cheb_tp}
There is a wide range of numerical methods which basically allow to map a sequence of SpMV steps to an MPK, such as s-step Krylov methods, exponential time differencing, polynomial preconditioning, or eigenvalue computations.
Here we focus on the solution of time-dependent (dynamic) partial differential equations (PDEs) using an exponential time evolution operator $U(\Delta t)$~\cite{Suhov} .
We choose Chebyshev polynomials to approximate the exponential, i.e., 
%
$U(\Delta t) = \sum_{k=0}^{M} c_k(\Delta t) T_k(A)$, where $c_k$ are  the coefficients (which depend on the time step $\Delta t$), $M$ is the number of Chebyshev  moments, $A$ is a sparse matrix derived from the underlying PDE, and $T_k(A)$ are Chebyshev polynomials of order $k$.
In our applications, Chebyshev polynomials of the first kind are used; therefore, $T_k(A)$ is defined using the following recurrence relation:
\begin{align}
& T_0(A) = I \text{, } T_1(A) = A \text{, }\nonumber \\
& T_{k+1}(A) = 2AT_{k}(A) - T_{k-1}(A) \eos \label{eq:cheb_rec}
\end{align}
The time evolution can then be computed by applying the operator $U(\Delta t)$ to 
the current state vector $x_t$ to obtain the next state vector $x_{t+\Delta t}$, i.e., 
\begin{equation}
x_{t+\Delta t} = U(\Delta t)x_t \approx \sum_{k=0}^{M} c_k(\Delta t) v_k \cma \label{eq:cheb_tp_time_evolution}
\end{equation}
 where $v_k = T_k(A)x_t$. 
The polynomial matrices $T_k(A)$ need not be stored explicitly as the $v_{k+1}$ can be determined from previous $v_k$ exploiting the recurrence relation~\eqref{eq:cheb_rec}, i.e., 
\begin{align}
v_{k+1} = 2Av_{k} - v_{k-1} \eos \label{eq:chebvec_rec}
\end{align}
Thus, the computation of $U(\Delta t)x_t$ can be implemented as a sequence of $M$ SpMVs ($Av_{k}$) with appropriate scaling factors. 
This step, which propagates the system by one time step, is the compute time hotspot of Cheb-TP applications.
Current  implementations of Cheb-TP (see, e.g.,~\cite{FEHSKE2009}, \cite{Suhov}, \cite{SCHAEFER2017}, \cite{recent_chebtp}) perform successive back-to-back SpMVs similar to our baseline MPK.
However, the relation \eqref{eq:chebvec_rec} allows to apply the RACE MPK in the computation of $U(\Delta t)x_t$ and to perform level-based cache blocking across $p_m$ successive SpMVs. 
As Cheb-TP typically uses high values of $M$ (few 100s--1000s), choosing $p_m = M$ is not advisable (see Sec.~\ref{sec:power_influence}).
Hence, we split the $M$ SpMVs into $M/p_m$ batches and perform $p_m$ successive SpMVs via the RACE MPK within each batch.
\begin{table}[tbp]
	\centering
	\caption{Details of the matrices used for the Cheb-TP application.\label{tab:matrices_app}}
	\begin{center}
		\resizebox{\linewidth}{!}{%
			\begin{tabular}{|l|l|S[table-format=7.0, table-space-text-pre=(, table-space-text-post=)]|S[table-format=8.0, table-space-text-pre=(, table-space-text-post=)]|S[round-mode=places,round-precision=2]|}
\toprule
{Matrix name} &  {Properties} & {\NR} & {\NNZ} & {\NNZR} \\
\midrule
{order-2} & grid size=$160^3$, spatial order=2 & 4096000 & 28518400 & 6.962500 \\
{order-4} & grid size=$160^3$, spatial  order=4 & 4096000 & 52787200 & 12.887500\\
{order-6} & grid size=$160^3$, spatial  order=6  & 4096000 & 76902400 & 18.775000\\
\midrule
{Fermion} & number of sites=24, fermions=12 & 2704156 & 32449872 & 12 \\
{Graphene} & sheet size=$2048^2$ & 4194304 & 54480904 & 12.98925972\\
{Anderson} & lattice sites=$200^3$ & 8000000 & 55760000 & 6.97\\
\bottomrule
\end{tabular}

		}
	\end{center}
\end{table}

\begin{figure*}[tbp]
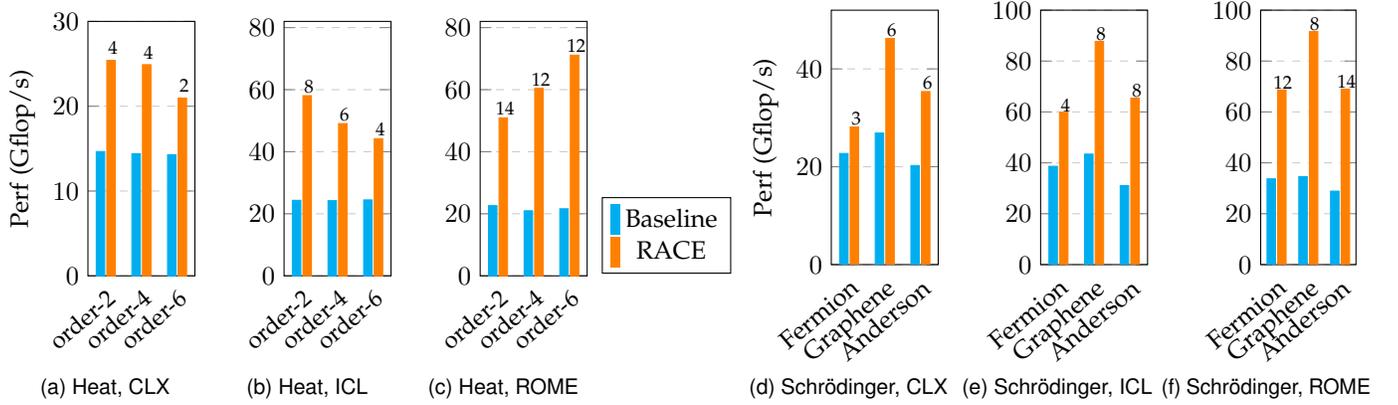

	\centering
	\subfloat[Heat, CLX]{%
	\input{plots/tikz/cheb_tp/plot_casclakesp2_heat.tex}%

		\label{fig:cheb_tp:heat:clx}}
	\subfloat[Heat, ICL]{%
	\input{plots/tikz/cheb_tp/plot_horeka_heat.tex}%

		\label{fig:cheb_tp:heat:icl}}
	\subfloat[Heat, ROME]{%
	\input{plots/tikz/cheb_tp/plot_tg094_heat.tex}%

		\label{fig:cheb_tp:heat:rome}}
	%
	\input{plots/tikz/cheb_tp/cheb_tp_legend.tex}%

	\subfloat[Schr\"{o}dinger, CLX]{%
	\input{plots/tikz/cheb_tp/plot_casclakesp2_schroedinger.tex}%

		\label{fig:cheb_tp:schroedinger:clx}}
	\subfloat[Schr\"{o}dinger, ICL]{%
	\input{plots/tikz/cheb_tp/plot_horeka_schroedinger.tex}%

		\label{fig:cheb_tp:schroedinger:icl}}
	\subfloat[Schr\"{o}dinger, ROME]{%
	\input{plots/tikz/cheb_tp/plot_tg094_schroedinger.tex}%

		\label{fig:cheb_tp:schroedinger:rome}}
	\caption{Cheb-TP performance using the baseline and the RACE MPK on the three architectures.
		(a--c) Three cases of 3D heat equation with various spatial discretization order.
		(d--f) Schr\"{o}dinger equation with three different Hamiltonian matrices $H$ (see \eqref{eq:schroedinger}) derived from various physics applications. 
		In all cases the tuned $p_m$ values of RACE are  annotated in the plots.
		Note the difference in $y$-axis scaling for CLX data compared to ICL and ROME. 
			\label{fig:cheb_tp}}
\end{figure*}
To demonstrate the performance potential of RACE MPK, we perform time propagation for PDEs underlying the parabolic heat equation 
\begin{equation}
	\frac{\partial \psi(x,t)}{\partial t} = a\Delta \psi(x,t) \label{eq:heat}
\end{equation}
and the Schr\"{o}dinger equation
\begin{equation}
	i\hbar\frac{\partial \psi(x,t)}{\partial t} = H\psi(x,t)\eos \label{eq:schroedinger}
\end{equation}
choosing three application scenarios (i.e., representative sparse matrices) for both.
These scenarios have been selected such that the system matrices $A$ are sparse and real. 
For the heat equation~\eqref{eq:heat}, the coefficients $c_k$ depend on the modified Bessel functions and are similar to the ones found in~\cite{tal_ezer_parabolic}.
This means the code balance of the entire algorithm implemented using the baseline MPK approach
would be similar to that of SpMV as discussed in Sec.~\ref{sec:baseline_MPK}.
For the Schr\"{o}dinger equation~\eqref{eq:schroedinger}, the coefficients are complex numbers and depend on the Bessel functions~\cite{FEHSKE2009}.
The complex coefficients result in complex vectors $v_k$, and thus SpMVs between a real matrix and complex vectors have to be performed.
This results in nearly twice as many flops due to complex arithmetic, but the data traffic remains almost the same as the matrix has real entries.
Therefore, the code balance is reduced by a factor of roughly two in case of the Schr\"{o}dinger equation.

The value of $M$ in Cheb-TP method has been chosen such that the Bessel coefficients have a value lower than $10^{-4}$ (see the cut-off strategy in~\cite{FEHSKE2009}).
For RACE MPK, the tuning space is similar to the one used in Sec.~\ref{sec:mpk:exp_setup} above.
Figure~\ref{fig:cheb_tp} shows the  performance in \GFS\ of Cheb-TP  using the baseline and RACE MPK approaches. 
Note that the baseline performance of Cheb-TP on the Schr\"{o}dinger equation is almost twice that of the heat equation due to 2$\times$ lower code balance.
In general it can be seen  that level-blocked approach of RACE outperforms the baseline method of back-to-back SpMVs. 
In line with our results in Sec.~\ref{sec:perf_eval}, the highest speedup is attained on ROME followed by ICL and CLX.

For the heat equation we use a three-dimensional grid of size $160^3$ and discretize the spatial derivatives (i.e., $\Delta \psi$ in \eqref{eq:heat}) using finite differences of second, fourth, and sixth order (order-2 \ldots\ order-6 in Fig.~\ref{fig:cheb_tp}).
Of course, higher-order discretization leads to an increase 
in the number of neighbors, i.e., \NNZR\ of the sparse matrix (see Table~\ref{tab:matrices_app} for details).
This in turn results in bulkier levels and therefore 
should lead to a lower optimal $p_m$ value due to~\eqref{eq:lc} and the recursion overhead.
This results in lower performance of RACE MPK with increasing discretization order, as can be seen in Fig.~\ref{fig:cheb_tp:heat:clx} and  Fig.~\ref{fig:cheb_tp:heat:icl} for CLX and ICL.
However, the large cache size of ROME allows to maintain high $p_m$ values even for large 
orders and thus performance increases with discretization order as seen in Fig.~\ref{fig:cheb_tp:heat:rome}.
This is due to the increase in \NNZR\  causing a relative reduction in  vector traffic contributions and allowing for efficient SIMD vectorization along the inner loop over \NNZR.
Overall, for the heat equation the RACE MPK attains an average speedup of 
2.8$\times$, 2.1$\times$, and 1.6$\times$ over the baseline approach on ROME, ICL, and CLX, respectively.

For the Schr\"odinger equation, we test our approach using three Hamiltonian 
matrices ($H$ in~\eqref{eq:schroedinger}) from quantum physics applications.
The Fermion and Anderson matrices model noninteracting many-particle quantum systems and 
the motion of a quantum-mechanical particle in a disordered solid~\cite{ScaMaC}, respectively.
Both matrices were generated using the ScaMaC~\cite{ScaMaC} library.
The Graphene matrix arises from modeling graphene tubes and ribbons~\cite{graphene}
and can be generated using the ESSEX-Physics library~\cite{essex_physics}.
The details of these matrices can be found in Table~\ref{tab:matrices_app}.
They have varying sparsity patterns resulting in widely different level structures and indirect access patterns, which results in a wide performance range (compared to the heat equation) of the baseline method for the three matrices (see Figs.~\ref{fig:cheb_tp:schroedinger:clx}, \ref{fig:cheb_tp:schroedinger:icl} and \ref{fig:cheb_tp:schroedinger:rome}).
Similar to the heat equation, when compared to the baseline approach the RACE MPK achieves an average speedup of 2.3$\times$, 1.9$\times$, and 1.6$\times$ on ROME, ICL, and CLX, respectively.



\vfill

\section{Conclusion and outlook}
\label{sec:conclusion}
In this paper we have developed a level-based blocking algorithm (RACE MPK) to increase the 
performance of sparse matrix-power-vector kernels (MPK). 
The RACE algorithm uses levels, generated by breadth-first search, to increase temporal access locality for the matrix entries by reusing them for successive power computations.
Various hardware-oriented algorithmic optimization strategies such as level grouping, 
point-to-point synchronization, and recursive application of the level-blocking scheme are introduced to further improve the performance of RACE MPK.
A thorough performance analysis on a representative set of 36 matrices 
shows that RACE MPK outperforms a standard MPK implementation
by an average factor of 2$\times$ and 3.5$\times$ on modern Intel and AMD CPUs.
Finally, applying RACE MPK to Chebyshev time propagation we demonstrate that similar speedups are  achievable for real-world applications.

The MPK finds its use in a large variety of applications, especially in
the field of communication-avoiding algorithms~\cite{hoemmen_thesis}, polynomial preconditioning~\cite{Jenifer_poly_precon}, and exponential time integration~\cite{Suhov}.
Future work includes integrating RACE MPK into
Krylov solvers and preconditioners from the Trilinos~\cite{trilinos-website} framework.

%


%




\ifCLASSOPTIONcompsoc
  \section*{Acknowledgments}
\else
  \section*{Acknowledgment}
\fi
This work was partially supported by NHR@FAU, which is funded by the State of Bavaria and  by
the Federal Ministry of Education and Research.
The authors would also like to thank NHR@KIT 
for providing access to the HoreKa supercomputer (ICL system), which is funded by the 
Ministry of Science, Research and the Arts Baden-Württemberg and by
the Federal Ministry of Education and Research.
The authors thank Kengo Nakajima for helpful discussions within the JHPCN project ``Innovative Multigrid Methods II.''

\newpage

\ifCLASSOPTIONcaptionsoff
  \newpage
\fi



%


\bibliographystyle{IEEEtran}
\bibliography{pub.bib}

\begin{thebibliography}{10}
\providecommand{\url}[1]{#1}
\csname url@samestyle\endcsname
\providecommand{\newblock}{\relax}
\providecommand{\bibinfo}[2]{#2}
\providecommand{\BIBentrySTDinterwordspacing}{\spaceskip=0pt\relax}
\providecommand{\BIBentryALTinterwordstretchfactor}{4}
\providecommand{\BIBentryALTinterwordspacing}{\spaceskip=\fontdimen2\font plus
\BIBentryALTinterwordstretchfactor\fontdimen3\font minus
  \fontdimen4\font\relax}
\providecommand{\BIBforeignlanguage}[2]{{%
\expandafter\ifx\csname l@#1\endcsname\relax
\typeout{** WARNING: IEEEtran.bst: No hyphenation pattern has been}%
\typeout{** loaded for the language `#1'. Using the pattern for}%
\typeout{** the default language instead.}%
\else
\language=\csname l@#1\endcsname
\fi
#2}}
\providecommand{\BIBdecl}{\relax}
\BIBdecl

\bibitem{GPU_review}
\BIBentryALTinterwordspacing
S.~Filippone, V.~Cardellini, D.~Barbieri, and A.~Fanfarillo, ``Sparse
  matrix-vector multiplication on {GPGPUs},'' \emph{ACM Trans. Math. Softw.},
  vol.~43, no.~4, jan 2017. [Online]. Available:
  \url{https://doi.org/10.1145/3017994}
\BIBentrySTDinterwordspacing

\bibitem{barrett_formats}
\BIBentryALTinterwordspacing
R.~Barrett, M.~Berry, T.~F. Chan, J.~Demmel, J.~Donato, J.~Dongarra,
  V.~Eijkhout, R.~Pozo, C.~Romine, and H.~van~der Vorst, \emph{Templates for
  the Solution of Linear Systems: Building Blocks for Iterative Methods}.\hskip
  1em plus 0.5em minus 0.4em\relax Society for Industrial and Applied
  Mathematics, 1994. [Online]. Available:
  \url{https://epubs.siam.org/doi/abs/10.1137/1.9781611971538}
\BIBentrySTDinterwordspacing

\bibitem{Vuduc2003:thesis}
R.~W. Vuduc, ``Automatic performance tuning of sparse matrix kernels,'' Ph.D.
  dissertation, University of California, Berkeley, December 2003.

\bibitem{CSB}
\BIBentryALTinterwordspacing
A.~Bulu\c{c}, J.~T. Fineman, M.~Frigo, J.~R. Gilbert, and C.~E. Leiserson,
  ``Parallel sparse matrix-vector and matrix-transpose-vector multiplication
  using compressed sparse blocks,'' in \emph{Proceedings of the Twenty-First
  Annual Symposium on Parallelism in Algorithms and Architectures}, ser. SPAA
  '09.\hskip 1em plus 0.5em minus 0.4em\relax New York, NY, USA: Association
  for Computing Machinery, 2009, p. 233–244. [Online]. Available:
  \url{https://doi.org/10.1145/1583991.1584053}
\BIBentrySTDinterwordspacing

\bibitem{KreutzerSELLC14}
\BIBentryALTinterwordspacing
M.~Kreutzer, G.~Hager, G.~Wellein, H.~Fehske, and A.~R. Bishop, ``A unified
  sparse matrix data format for efficient general sparse matrix-vector
  multiplication on modern processors with wide {SIMD} units,'' \emph{SIAM
  Journal on Scientific Computing}, vol.~36, no.~5, pp. C401--C423, 2014.
  [Online]. Available: \url{https://doi.org/10.1137/130930352}
\BIBentrySTDinterwordspacing

\bibitem{SpMV_reordering}
\BIBentryALTinterwordspacing
L.~Oliker, X.~Li, P.~Husbands, and R.~Biswas, ``Effects of ordering strategies
  and programming paradigms on sparse matrix computations,'' \emph{SIAM
  Review}, vol.~44, no.~3, pp. 373--393, 2002. [Online]. Available:
  \url{https://doi.org/10.1137/S00361445003820}
\BIBentrySTDinterwordspacing

\bibitem{BEBOP2006}
R.~Nishtala, R.~W. Vuduc, J.~W. Demmel, and K.~Yelick, ``{When Cache Blocking
  Sparse Matrix Vector Multiply Works and Why},'' \emph{Applicable Algebra in
  Engineering, Communication and Computing}, vol.~18, no.~3, pp. 297--311,
  2007.

\bibitem{balaprakash2018autotuning}
P.~Balaprakash, J.~Dongarra, T.~Gamblin, M.~Hall, J.~K. Hollingsworth,
  B.~Norris, and R.~Vuduc, ``Autotuning in high-performance computing
  applications,'' \emph{Proceedings of the IEEE}, vol. 106, no.~11, pp.
  2068--2083, 2018.

\bibitem{hong2019adaptive}
C.~Hong, A.~Sukumaran-Rajam, I.~Nisa, K.~Singh, and P.~Sadayappan, ``Adaptive
  sparse tiling for sparse matrix multiplication,'' in \emph{Proceedings of the
  24th Symposium on Principles and Practice of Parallel Programming}, 2019, pp.
  300--314.

\bibitem{RACE}
\BIBentryALTinterwordspacing
C.~Alappat, A.~Basermann, A.~R. Bishop, H.~Fehske, G.~Hager, O.~Schenk,
  J.~Thies, and G.~Wellein, ``A recursive algebraic coloring technique for
  hardware-efficient symmetric sparse matrix-vector multiplication,'' \emph{ACM
  Trans. Parallel Comput.}, vol.~7, no.~3, Jun. 2020. [Online]. Available:
  \url{https://doi.org/10.1145/3399732}
\BIBentrySTDinterwordspacing

\bibitem{demmel_MPK_tr}
\BIBentryALTinterwordspacing
J.~Demmel, M.~F. Hoemmen, M.~Mohiyuddin, and K.~A. Yelick, ``Avoiding
  communication in computing {K}rylov subspaces,'' EECS Department, University
  of California, Berkeley, Tech. Rep. UCB/EECS-2007-123, Oct 2007. [Online].
  Available:
  \url{http://www2.eecs.berkeley.edu/Pubs/TechRpts/2007/EECS-2007-123.html}
\BIBentrySTDinterwordspacing

\bibitem{hoemmen_thesis}
M.~Hoemmen, ``Communication-avoiding krylov subspace methods,'' Ph.D.
  dissertation, USA, 2010, aAI3413388.

\bibitem{Carson:EECS-2015-179}
\BIBentryALTinterwordspacing
E.~Carson, ``Communication-avoiding krylov subspace methods in theory and
  practice,'' Ph.D. dissertation, EECS Department, University of California,
  Berkeley, Aug 2015. [Online]. Available:
  \url{http://www2.eecs.berkeley.edu/Pubs/TechRpts/2015/EECS-2015-179.html}
\BIBentrySTDinterwordspacing

\bibitem{7914608}
J.~Dongarra, S.~Tomov, P.~Luszczek, J.~Kurzak, M.~Gates, I.~Yamazaki, H.~Anzt,
  A.~Haidar, and A.~Abdelfattah, ``With extreme computing, the rules have
  changed,'' \emph{Computing in Science Engineering}, vol.~19, no.~3, pp.
  52--62, 2017.

\bibitem{mohiyudeen}
\BIBentryALTinterwordspacing
M.~Mohiyuddin, M.~Hoemmen, J.~Demmel, and K.~Yelick, ``Minimizing communication
  in sparse matrix solvers,'' in \emph{Proceedings of the Conference on High
  Performance Computing Networking, Storage and Analysis}, ser. SC '09.\hskip
  1em plus 0.5em minus 0.4em\relax New York, NY, USA: Association for Computing
  Machinery, 2009. [Online]. Available:
  \url{https://doi.org/10.1145/1654059.1654096}
\BIBentrySTDinterwordspacing

\bibitem{PITCH_temporal_overview}
\BIBentryALTinterwordspacing
T.~Muranushi and J.~Makino, ``Optimal temporal blocking for stencil
  computation,'' \emph{Procedia Computer Science}, vol.~51, pp. 1303--1312,
  2015, international Conference On Computational Science, ICCS 2015. [Online].
  Available:
  \url{https://www.sciencedirect.com/science/article/pii/S1877050915011230}
\BIBentrySTDinterwordspacing

\bibitem{schreiber_multicore}
\BIBentryALTinterwordspacing
D.~Huber, M.~Schreiber, and M.~Schulz, ``Graph-based multi-core higher-order
  time integration of linear autonomous partial differential equations,''
  \emph{Journal of Computational Science}, vol.~53, p. 101349, 2021. [Online].
  Available:
  \url{https://www.sciencedirect.com/science/article/pii/S1877750321000466}
\BIBentrySTDinterwordspacing

\bibitem{Ichi_GMRES_paper}
I.~Yamazaki, H.~Anzt, S.~Tomov, M.~Hoemmen, and J.~Dongarra, ``Improving the
  performance of {CA}-{GMRES} on multicores with multiple {GPUs},'' in
  \emph{2014 IEEE 28th International Parallel and Distributed Processing
  Symposium}, 2014, pp. 382--391.

\bibitem{Ichi_precon}
I.~Yamazaki, S.~Rajamanickam, E.~G. Boman, M.~Hoemmen, M.~A. Heroux, and
  S.~Tomov, ``Domain decomposition preconditioners for communication-avoiding
  {K}rylov methods on a hybrid {CPU}/{GPU} cluster,'' in \emph{SC '14:
  Proceedings of the International Conference for High Performance Computing,
  Networking, Storage and Analysis}, 2014, pp. 933--944.

\bibitem{diamond_MPK}
\BIBentryALTinterwordspacing
E.~Vatai, U.~Singhal, and R.~Suda, ``Diamond matrix powers kernels,'' in
  \emph{Proceedings of the International Conference on High Performance
  Computing in Asia-Pacific Region}, ser. HPCAsia2020.\hskip 1em plus 0.5em
  minus 0.4em\relax New York, NY, USA: Association for Computing Machinery,
  2020, p. 102–113. [Online]. Available:
  \url{https://doi.org/10.1145/3368474.3368494}
\BIBentrySTDinterwordspacing

\bibitem{doi:10.1137/070693199}
\BIBentryALTinterwordspacing
K.~Datta, S.~Kamil, S.~Williams, L.~Oliker, J.~Shalf, and K.~Yelick,
  ``Optimization and performance modeling of stencil computations on modern
  microprocessors,'' \emph{SIAM Review}, vol.~51, no.~1, pp. 129--159, 2009.
  [Online]. Available: \url{https://doi.org/10.1137/070693199}
\BIBentrySTDinterwordspacing

\bibitem{6012879}
M.~Christen, O.~Schenk, and H.~Burkhart, ``{PATUS}: A code generation and
  autotuning framework for parallel iterative stencil computations on modern
  microarchitectures,'' in \emph{2011 IEEE International Parallel Distributed
  Processing Symposium}, 2011, pp. 676--687.

\bibitem{9355259}
H.~Wang and A.~Chandramowlishwaran, ``Pencil: A pipelined algorithm for
  distributed stencils,'' in \emph{SC20: International Conference for High
  Performance Computing, Networking, Storage and Analysis}, 2020, pp. 1--16.

\bibitem{10.1145/3293883.3295712}
\BIBentryALTinterwordspacing
C.~Hong, A.~Sukumaran-Rajam, I.~Nisa, K.~Singh, and P.~Sadayappan, ``Adaptive
  sparse tiling for sparse matrix multiplication,'' in \emph{Proceedings of the
  24th Symposium on Principles and Practice of Parallel Programming}, ser.
  PPoPP '19.\hskip 1em plus 0.5em minus 0.4em\relax New York, NY, USA:
  Association for Computing Machinery, 2019, p. 300–314. [Online]. Available:
  \url{https://doi.org/10.1145/3293883.3295712}
\BIBentrySTDinterwordspacing

\bibitem{top500_list}
\BIBentryALTinterwordspacing
``{Top 500: June 2021 list}.'' [Online]. Available:
  \url{https://www.top500.org/lists/top500/2021/06/}
\BIBentrySTDinterwordspacing

\bibitem{Understanding}
C.~L. Alappat, J.~Hofmann, G.~Hager, H.~Fehske, A.~R. Bishop, and G.~Wellein,
  ``Understanding {HPC} benchmark performance on {Intel} {Broadwell} and
  {Cascade} {Lake} processors,'' in \emph{High Performance Computing},
  P.~Sadayappan, B.~L. Chamberlain, G.~Juckeland, and H.~Ltaief, Eds.\hskip 1em
  plus 0.5em minus 0.4em\relax Cham: Springer International Publishing, 2020,
  pp. 412--433.

\bibitem{UOF}
\BIBentryALTinterwordspacing
T.~A. Davis and Y.~Hu, ``{The University of Florida Sparse Matrix
  Collection},'' \emph{ACM Trans. Math. Softw.}, vol.~38, no.~1, pp. 1:1--1:25,
  Dec. 2011. [Online]. Available:
  \url{http://doi.acm.org/10.1145/2049662.2049663}
\BIBentrySTDinterwordspacing

\bibitem{HPCG}
\BIBentryALTinterwordspacing
M.~A. Heroux and J.~Dongarra, ``Toward a new metric for ranking high
  performance computing systems.'' 6 2013, {HPCG} website:
  \url{https://www.hpcg-benchmark.org}. [Online]. Available:
  \url{https://www.osti.gov/biblio/1089988}
\BIBentrySTDinterwordspacing

\bibitem{CRS_saad}
Y.~Saad, ``{SPARSKIT}: a basic tool kit for sparse matrix computations,''
  {Research Institute for Advanced Computer Science}, Tech. Rep., 1990.

\bibitem{Gropp99towardsrealistic}
W.~D. Gropp, D.~K. Kaushik, D.~E. Keyes, and B.~F. Smith, ``Towards realistic
  performance bounds for implicit {CFD} codes,'' in \emph{Proceedings of
  Parallel CFD’99}.\hskip 1em plus 0.5em minus 0.4em\relax Elsevier, 1999,
  pp. 233--240.

\bibitem{Frigo}
\BIBentryALTinterwordspacing
M.~Frigo and V.~Strumpen, ``The memory behavior of cache oblivious stencil
  computations,'' \emph{The Journal of Supercomputing}, vol.~39, no.~2, pp.
  93--112, 2007. [Online]. Available:
  \url{https://doi.org/10.1007/s11227-007-0111-y}
\BIBentrySTDinterwordspacing

\bibitem{diamond_blocking}
\BIBentryALTinterwordspacing
T.~Malas, G.~Hager, H.~Ltaief, H.~Stengel, G.~Wellein, and D.~Keyes,
  ``Multicore-optimized wavefront diamond blocking for optimizing stencil
  updates,'' \emph{SIAM Journal on Scientific Computing}, vol.~37, no.~4, pp.
  C439--C464, 2015. [Online]. Available:
  \url{https://doi.org/10.1137/140991133}
\BIBentrySTDinterwordspacing

\bibitem{BFS}
C.~Y. Lee, ``An algorithm for path connections and its applications,''
  \emph{IRE Transactions on Electronic Computers}, vol. EC-10, no.~3, pp.
  346--365, Sept 1961.

\bibitem{RCM_SpMP}
K.~I. Karantasis, A.~Lenharth, D.~Nguyen, M.~J. Garzarán, and K.~Pingali,
  ``Parallelization of reordering algorithms for bandwidth and wavefront
  reduction,'' in \emph{SC '14: Proceedings of the International Conference for
  High Performance Computing, Networking, Storage and Analysis}, 2014, pp.
  921--932.

\bibitem{Malas2017}
\BIBentryALTinterwordspacing
T.~M. Malas, G.~Hager, H.~Ltaief, and D.~E. Keyes, ``Multidimensional intratile
  parallelization for memory-starved stencil computations,'' \emph{ACM Trans.
  Parallel Comput.}, vol.~4, no.~3, pp. 12:1--12:32, Dec. 2017. [Online].
  Available: \url{http://doi.acm.org/10.1145/3155290}
\BIBentrySTDinterwordspacing

\bibitem{perf_pattern}
T.~R{\"o}hl, J.~Eitzinger, G.~Hager, and G.~Wellein, ``Validation of hardware
  events for successful performance pattern identification in high performance
  computing,'' in \emph{Tools for High Performance Computing 2015},
  A.~Kn{\"u}pfer, T.~Hilbrich, C.~Niethammer, J.~Gracia, W.~E. Nagel, and M.~M.
  Resch, Eds.\hskip 1em plus 0.5em minus 0.4em\relax Cham: Springer
  International Publishing, 2016, pp. 17--28.

\bibitem{RACE-git}
\BIBentryALTinterwordspacing
C.~Alappat, \emph{Recursive Algebraic Coloring Engine library}, 2019 (acccessed
  May 2, 2022). [Online]. Available: \url{https://github.com/RRZE-HPC/RACE}
\BIBentrySTDinterwordspacing

\bibitem{SpMP}
\BIBentryALTinterwordspacing
{SpMP Development Team}, ``Sparse matrix pre-processing library.'' [Online].
  Available: \url{https://github.com/IntelLabs/SpMP}
\BIBentrySTDinterwordspacing

\bibitem{likwid_L3_problem}
\BIBentryALTinterwordspacing
``{L2 L3 MEM traffic on Intel Skylake SP CascadeLake SP}.'' [Online].
  Available:
  \url{https://github.com/RRZE-HPC/likwid/wiki/L2-L3-MEM-traffic-on-Intel-Skylake-SP-CascadeLake-SP}
\BIBentrySTDinterwordspacing

\bibitem{BFS_Scott}
S.~Beamer, K.~Asanovic, and D.~Patterson, ``Direction-optimizing breadth-first
  search,'' in \emph{SC '12: Proceedings of the International Conference on
  High Performance Computing, Networking, Storage and Analysis}, 2012, pp.
  1--10.

\bibitem{Suhov}
\BIBentryALTinterwordspacing
A.~Y. Suhov, ``An accurate polynomial approximation of exponential
  integrators,'' \emph{Journal of Scientific Computing}, vol.~60, no.~3, pp.
  684--698, 2014. [Online]. Available:
  \url{https://doi.org/10.1007/s10915-013-9813-x}
\BIBentrySTDinterwordspacing

\bibitem{FEHSKE2009}
\BIBentryALTinterwordspacing
H.~Fehske, J.~Schleede, G.~Schubert, G.~Wellein, V.~S. Filinov, and A.~R.
  Bishop, ``Numerical approaches to time evolution of complex quantum
  systems,'' \emph{Physics Letters A}, vol. 373, no.~25, pp. 2182--2188, 2009.
  [Online]. Available:
  \url{https://www.sciencedirect.com/science/article/pii/S0375960109004927}
\BIBentrySTDinterwordspacing

\bibitem{SCHAEFER2017}
\BIBentryALTinterwordspacing
I.~Schaefer, H.~Tal-Ezer, and R.~Kosloff, ``Semi-global approach for
  propagation of the time-dependent {S}chr{\"o}dinger equation for
  time-dependent and nonlinear problems,'' \emph{Journal of Computational
  Physics}, vol. 343, pp. 368--413, 2017. [Online]. Available:
  \url{https://www.sciencedirect.com/science/article/pii/S0021999117302887}
\BIBentrySTDinterwordspacing

\bibitem{recent_chebtp}
T.~Löthman, C.~Triola, J.~Cayao, and A.~M. Black-Schaffer, ``Efficient
  numerical method for evaluating normal and anomalous time-dependent
  equilibrium {Green}'s functions in inhomogeneous systems,'' 2020.

\bibitem{tal_ezer_parabolic}
\BIBentryALTinterwordspacing
H.~Tal-Ezer, ``Spectral methods in time for parabolic problems,'' \emph{SIAM
  Journal on Numerical Analysis}, vol.~26, no.~1, pp. 1--11, 1989. [Online].
  Available: \url{https://doi.org/10.1137/0726001}
\BIBentrySTDinterwordspacing

\bibitem{ScaMaC}
{Andreas Alvermann}, ``{ScaMaC}: The scalable matrix collection,''
  \url{https://bitbucket.org/essex/matrixcollection/}, 2019.

\bibitem{graphene}
\BIBentryALTinterwordspacing
A.~H. Castro~Neto, F.~Guinea, N.~M.~R. Peres, K.~S. Novoselov, and A.~K. Geim,
  ``The electronic properties of graphene,'' \emph{Rev. Mod. Phys.}, vol.~81,
  pp. 109--162, Jan 2009. [Online]. Available:
  \url{https://link.aps.org/doi/10.1103/RevModPhys.81.109}
\BIBentrySTDinterwordspacing

\bibitem{essex_physics}
{ESSEX project team}, ``{ESSEX-Physics},''
  \url{https://bitbucket.org/essex/physics/src/master/}.

\bibitem{Jenifer_poly_precon}
\BIBentryALTinterwordspacing
J.~A. Loe, H.~K. Thornquist, and E.~G. Boman, \emph{Polynomial Preconditioned
  {GMRES} in Trilinos: Practical Considerations for High-Performance
  Computing}, pp. 35--45. [Online]. Available:
  \url{https://epubs.siam.org/doi/abs/10.1137/1.9781611976137.4}
\BIBentrySTDinterwordspacing

\bibitem{trilinos-website}
\BIBentryALTinterwordspacing
\relax The {T}rilinos~{P}roject {T}eam, \emph{The {T}rilinos {P}roject
  {W}ebsite}, 2021 (acccessed Aug 6, 2021). [Online]. Available:
  \url{https://trilinos.github.io}
\BIBentrySTDinterwordspacing

\end{thebibliography}

%
\newpage

\begin{IEEEbiography}[{\includegraphics[width=1.05in,height=1.9in,clip,keepaspectratio]{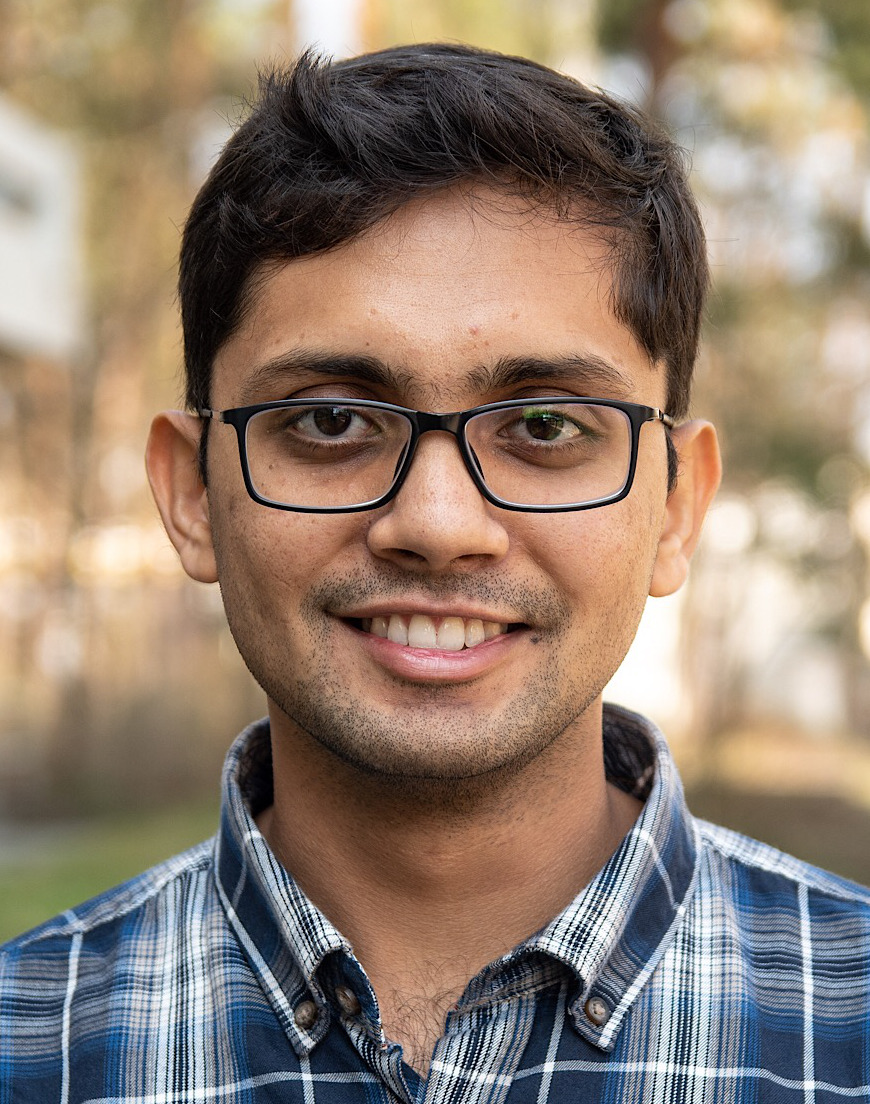}}]{Christie Alappat} received a 
Master’s Degree with honors from the Bavarian Graduate School of Computational Engineering at Friedrich-Alexander-Universität Erlangen-Nürnberg.
He is currently doing his Ph.D.\ under the guidance of Prof.\ Gerhard Wellein. 
His research interests include performance engineering, sparse matrix and graph algorithms, iterative linear solvers, and eigenvalue computations.
\end{IEEEbiography}

\vspace{-20em}

\begin{IEEEbiography}[{\includegraphics[width=1.05in,height=1.9in,clip,keepaspectratio]{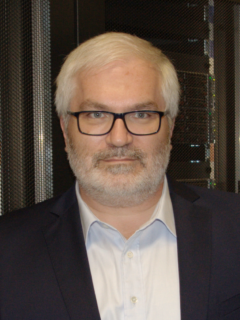}}]{Georg Hager} holds a doctorate (Ph.D.) and a Habilitation degree in Computational Physics from the University of Greifswald, Germany. He leads the Training \& Support Division at Erlangen National High Performance Computing Center (NHR@FAU) and is an associate lecturer at the Institute of Physics at the University of Greifswald. Recent research includes architecture-specific optimization strategies for current microprocessors, performance engineering of scientific codes on chip and system levels, and the modeling of out-of-lockstep behavior in large-scale parallel codes.
\end{IEEEbiography}

\vspace{-20em}

\begin{IEEEbiography}[{\includegraphics[width=1.05in,height=1.9in,clip,keepaspectratio]{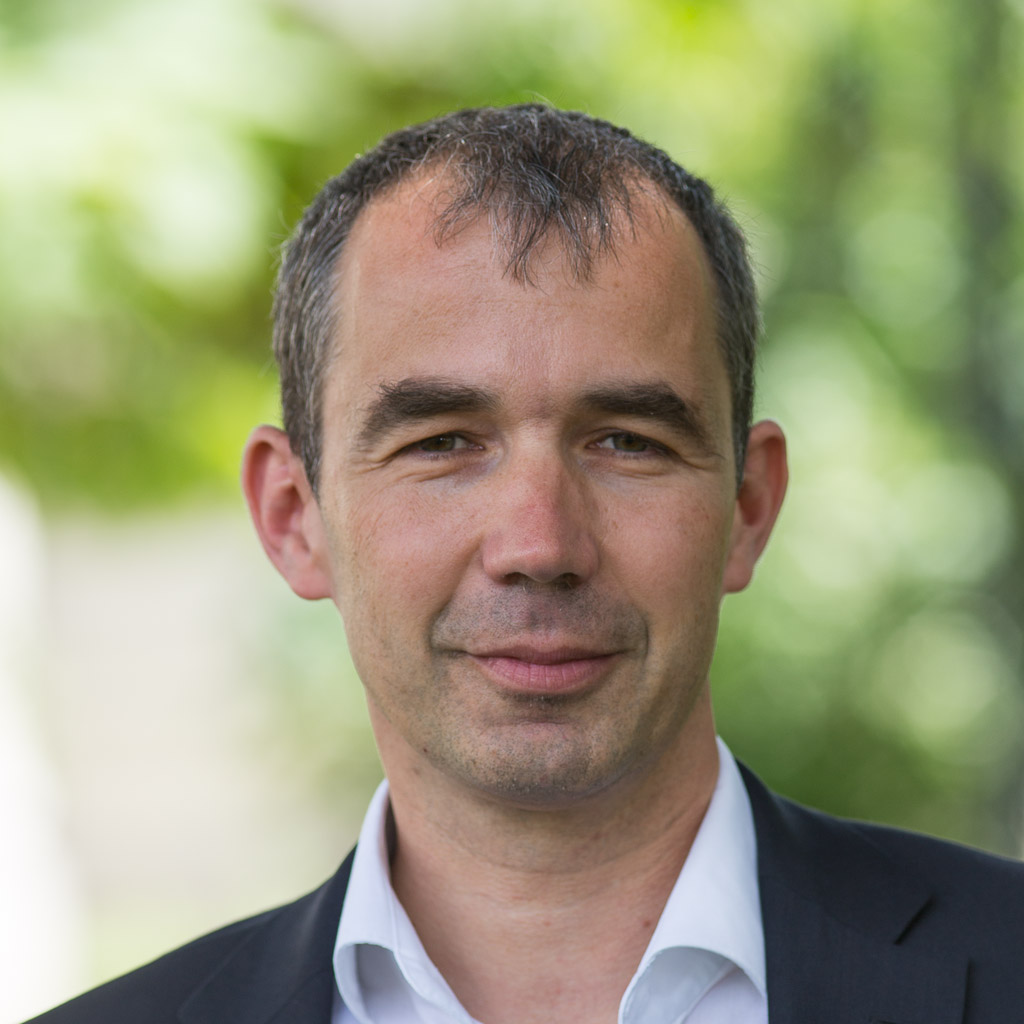}}]{Olaf Schenk}
	(M’02–SM’13) received the Diploma (M.Sc.) degree in mathematics from the University of Karlsruhe, Germany, and a doctorate (Ph.D.) degree in electrical engineering and information technology from the Swiss Federal Institute of Technology (ETH), Zurich, Switzerland. He is a Full Professor with the Institute of Computing within the Faculty of Informatics, Università della Svizzera italiana, Lugano, Switzerland, where he heads the Advanced Computing Laboratory. His research interests include extreme-scale simulations in computational algorithms, data science, application software, programming, and software tools.
\end{IEEEbiography}

\vspace{-20em}

\begin{IEEEbiography}[{\includegraphics[width=1.05in,height=1.9in,clip,keepaspectratio]{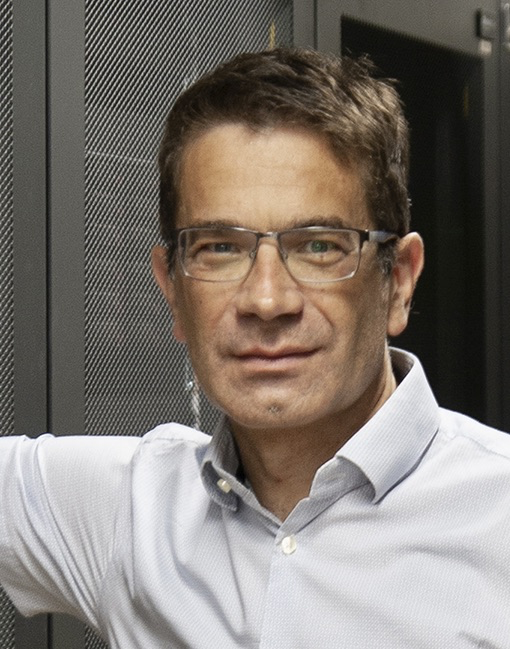}}]{Gerhard Wellein}
 received the Diploma (M.Sc.) degree and a doctorate (Ph.D.) degree in Physics from the University of Bayreuth, Germany. He is a Professor at the Department of Computer Science at Friedrich-Alexander-Universit\"at Erlangen-N\"urnberg and heads the Erlangen National Center for High-Performance Computing (NHR@FAU). His research interests focus on performance modeling and performance engineering, architecture-specific code optimization, and hardware-efficient building blocks for sparse linear algebra and stencil solvers.
\end{IEEEbiography}




\end{document}